\documentclass[12pt]{article}
\usepackage{amssymb}

\newtheorem{teo}{Theorem}
\newtheorem{prop}{Proposition}
\newtheorem{defi}{Definition}

\newtheorem{z}{Remark}
\newtheorem{lem}{Lemma}

\title{Uncertainty constants and quasispline wavelets\thanks{The work is supported by RFBR (grant 09-01-00162).}}
\date{ealebedeva2004@gmail.com}

\author{E.A.Lebedeva %\fnref{auth}
}
\begin{document}

\maketitle

\newcommand{\w}{\omega}
\newcommand{\ve}{\varepsilon}

\begin{abstract}
In 1996 Chui and Wang proved that the uncertainty constants of scaling and wavelet functions
tend to infinity as smoothness
of the wavelets grows for a broad class of wavelets such as Daubechies wavelets and spline  wavelets. We construct a  class of  new families of wavelets  (quasispline wavelets) whose uncertainty
constants tend to those of the Meyer wavelet function used in construction. 
\end{abstract}





\section{Introduction}
\label{intr}

One of the main advantages of wavelet systems  is the good time-frequency localization. The smoothness of wavelets is also a useful and desired property. So 
\begin{equation}
	\label{prob}
	\begin{array}{c}
	\mbox{\textit{to find wavelets preserving time-frequency}} \\
	\mbox{\textit{localization as smoothness grows}}\\
	\end{array}
\end{equation}
 is a  very attractive and interesting problem. In the sequel, by a wavelet we mean a function generating an orthonormal basis of $L_2(\mathbb{R})$ (see the definition in  section \ref{note}).  The measure of the time-frequency localization is  an uncertainty constant (see the definition in  section \ref{note}). So in problem (\ref{prob}) we are interested in uncertainty constants bounded with respect to a smoothness parameter.
 It is well known that the main classical families of wavelets contain wavelet functions with arbitrary large finite smoothness. Thus, one can investigate how a functional defined on a family of wavelets depends on the smoothness of the wavelets. Let the functional be the uncertainty constant.   
Unfortunately, the main classical families of wavelets lose  the time-frequency localization as the smoothness of chosen wavelet function grows.
More precisely, Chui and Wang \cite{ChW} show that the uncertainty constants of scaling and wavelet functions
tend to infinity as  the smoothness
of the wavelets grows for a broad class of wavelets such as for example Daubechies wavelets and spline  wavelets. So Daubechies wavelets and spline  wavelets don't settle (\ref{prob}).

Later Chui,  Wang \cite{ChW2} and Goodman,  Lee \cite{GL} construct families  of nonorthogonal scaling functions and semi-orthogonal wavelet functions. These functions have  optimal uncertainty constants (in the sense of Heisenberg uncertainty principle) as the smoothness parameter tends to infinity.
% But for orthogonal scaling and wavelet functions even  boundedness of  uncertainty constant was a still open question. 
But nothing is said about orthogonal scaling and wavelet functions in \cite{ChW2} and \cite{GL}.

Trying to solve problem (\ref{prob}),
Novikov \cite{N1}, \cite{N2} 
%using a nontrivial bright idea solved the problem (\ref{prob}) for the autocorrelation %function of the scaling function, but not for the scaling function and the wavelet. He %
constructs a family of modified Daubechies wavelets. The wavelet functions  are compactly supported.
The squared module of the modified Daubechies mask is the
Bernstein polynomial. It  interpolates a piecewise linear function (instead of the characteristic function of an
interval as it is in the case of classical Daubechies).  The smoothness of 
the modified Daubechies wavelet grows as the order of the Bernstein polynomial increases.
The time-frequency localization of the auto\-correlation function which is constructed for the scaling function of this family is preserved with the growth of the smoothness. It is still an open question 
whether the modified Daubechies scaling  and wavelet functions preserve the time-frequency localization as 
the smoothness growing.
%The point is that the following unsolved problem appears. Whether one can prove that for %some type of convergence from $\lim_{n \to \infty}|m_n|^2 = P$ it follows that 
%$\lim_{n \to \infty} m_n = \sqrt{ P},$ where $m_n$ is a trigonometric polynomial to the %power $n.$

In  \cite{LSMZ}, the author constructs a new wavelet family solving  problem (\ref{prob}) for  scaling functions. New scaling functions decay exponentially 
and their Fourier transforms decay as $O(\omega^{-l}),$ like spline wavelets; the uncertainty constants of
the scaling functions are uniformly bounded with respect to the smoothness parameter $l$.
The construction is based on  de la Valle-Poussin means of a function closely connected with the Meyer mask. 

In the present paper, 
we construct a wide class of such  wavelets (see  Theorem \ref{mainteo}). The new wavelet function also decays exponentially at infinity and its Fourier transform decays as $O(\omega^{-l}),$ like spline wavelet; that is why it is named a quasispline wavelet function (see Definition \ref{qswdf}).  The construction is based on the linear method of summation satisfying some weak, easily done  conditions (see Theorem \ref{cond}). The wavelet system constructed in \cite{LSMZ} is an example of the  quasispline wavelets. 
It is proven that the quasispline wavelets solve  problem (\ref{prob}) for  scaling  and  wavelet functions. Moreover, since the uncertainty constant for the Meyer scaling and wavelet function is bounded, 
 a special property for the quasispline wavelets is proven. This property is stronger than the boundedness. Namely,
we establish the  convergence of the uncertainty constants defined for
the new scaling (wavelet) functions   to those of the Meyer scaling (wavelet) function 
used in construction as the smoothness parameter $l$ goes to infinity.
The latter result also means that the uncertainty constant is a continuous functional, where the variable of the functional is a non-orthogonal mask $m_l.$ We also 
estimate the rate of the convergence.
It is necessary to note that the construction of quasispline wavelets can be based not only on the Meyer mask but also on any smooth orthogonal mask $m$ such that $m(\w)=1$ if $|\w|<a$ and $m(\w)=0$ if $b<|\w|<\pi$ for some $\pi/3 \leq a<b<\pi.$

\section{Notations and auxiliary results}
\label{note}
%$\mathbb{N}$ , $\mathbb{Z}$, $\mathbb{R}$ ~--- ��������� ����
%�����������, ����� � �������������� ����� ��������������.

Denote by $[x]$ an integer part of a real number  $x.$
Denote by $C^{k}[a,b]$ a space of all  $k$ times continuously differentiable functions defined on the interval $[a,b]$ with norm  
$
\|f\|_{W^k_{\infty}}:=\sum_{j=0}^k\max_{x\in [a,b]}|f^{(j)}(x)|,
$
write
$C^0[a,b]=C[a,b]$ and $C[-\pi,\,\pi]=C.$
% For function $f \in C_{[a,b]}$ module of continuity is defined by
%$$
%\omega(h,f):=\sup_ {\begin{array}{c}
%  x,x+x_1\in[a,b] \\
%  0<|x_1|\leq h \\
%\end{array}}
%|f(x+x_1)-f(x)|
%$$
%and satisfies the property
%$\omega(h,f)\longrightarrow 0$ as $h \longrightarrow 0.$
%Denote by  
%$\chi_M$ the characteristic function of a set $M$;
%and by
%$\zeta$ the Riemann zeta-function:
%$\zeta(x):=\sum_{n=1}^{\infty}n^{-x}.$

We choose \texttt{the Fourier transform} and the reconstruction formula as
$$
\widehat{f}(\w):=\int_{\mathbb{R}}f(t)e^{-it\w}\,dt,\quad f(t):=\frac{1}{2\pi}
\int_{\mathbb{R}}\widehat{f}(\w)e^{it\w}\,d\w
$$
respectively.
For the Fourier series
$
f \sim \frac{a_0}{2}+ \sum_{n \in \mathbb{N}}a_n \cos n\w + b_n \sin n\w
$
the sequence $(\lambda_{n,k}),$ $k=1, \dots n,$ $n \in \mathbb{N}$
defines   \texttt{a linear method of summation} 
$$
u_n(f,\w):=\frac{a_0}{2}+ \sum_{k=1}^n\lambda_{n,k}(a_n \cos n\w + b_n \sin n\w)
=
\int_{-\pi}^{\pi} f(x) U_n(x,\w)\,dx,
$$
where 
$U_n(x,\w):=1/2+ \sum_{k=1}^{n}\lambda_{n,k} \cos k (x-\w)$
and terms
$$
a_n:=\frac{1}{\pi} \int_{\pi}^{\pi}f(\w) \cos\,n\w\, d\w, \quad
b_n:=\frac{1}{\pi} \int_{\pi}^{\pi}f(\w) \sin\,n\w\, d\w
$$
are the Fourier coefficients.
The following property holds true
\begin{equation}
\label{u'}
	u_n(f',\w)=(u_n(f,\w))'_{\w}.
\end{equation}

%������� �����-������� ������� $f$ ����� ���:
%$$
%V_N(f,\w):=\pi^{-1}\int_{-\pi}^{\pi}f(t+\w)P_N(t)\,dt,
%$$
% ���
%$P_N(t):=(\cos Nt-\cos 2Nt)N^{-1}(2\sin 0,5t)^{-2}$
%~--- ���� �����-�������.
%
%������� ������ ��� ������� $f$ ����� ���:
%$$
%\sigma_N(f,\w):=\pi^{-1}\int_{-\pi}^{\pi}f(t+\w)F_N(t)\,dt,
%$$
%���
%$F_N(t):=2(N+1)^{-1}(\sin 0,5(N+1)t)^2(2\sin 0,5t)^{-2}$
%~--- ���� ������.
%
%������� �����-������� � ������ ������� ������������:
%\begin{equation}
%V_N(f,\w)=2\sigma_{2N-1}(f,\w)-\sigma_{N-1}(f,\w).
%\label{vp_f}
%\end{equation}
%
%��� ������������� ������� �������� �����-������� ������� ��������� ���������.
%\begin{teo} \cite[c.123-127]{KS}
%����   $f\in C_{[-\pi,\pi]}$, ��
%$$
%\left\|V_{\left[\frac{N+1}{2}\right]}(f,\cdot)-f(\cdot)\right\|_{C[-\pi,\pi]}\leq
%K\omega\left(\frac1N,f\right),
%$$
%��� $K$ ~--- ���������� ��������� (��������,$K=44$).
%\end{teo}

A function $\psi$ is called \texttt{a wavelet function} if the functions $2^{j/2}\psi(2^j \cdot -k),$ $j,\,k \in \mathbb{Z}$ form an orthonormal basis of $L_2(\mathbb{R})$. 

Denote by $\theta(\w)$ some odd function equal to $\frac{\pi}{4}$ for $\w>\frac{\pi}{3}.$ Assume henceforth that $\theta(\w)$ is a non-decreasing twice continuously differentiable function. Denote by $\w_0$ some parameter that varies in the interval $\frac{\pi}{3}\leq \w_0 <\frac{\pi}{2}$  and put $\w_1:=\pi-\w_0.$ \texttt{A Meyer scaling function} $\varphi^M$     is defined by
$$
    \widehat{\varphi^M}(\w):=\left\{
\begin{array}{cc}
  1, & |\w|\leq 2\w_0, \\
  \cos\left(\frac{\pi}{4}+\theta\left(\frac{\pi}{3(\pi-2\w_0)}(|\w|-\pi)\right)\right),
   & 2\w_0<|\w|\leq 2\pi-2\w_0, \\
  0, & |\w|>2\pi-2\w_0. \\
\end{array}\right.
$$
\texttt{A Meyer mask} is a $2\pi$-periodic function defined on $[-\pi,\;\pi]$ as follows
$m^M(\w):=
\widehat{\varphi^M}(2\w).$
%\begin{equation}
%m^M(\w):=
%\widehat{\varphi^M}(2\w)=
%\left\{
%\begin{array}{cc}
%  1, & |\w|\leq \w_0, \\
%  \cos\left(\frac{\pi}{4}+\theta\left(\frac{\pi}{3(\pi-2\w_0)}(2|\w|-\pi)\right)\right),
%   & \w_0<|\w|\leq \pi-\w_0 ,\\
%  0, & \pi-\w_0<|\w|\leq \pi .\\
%\end{array}\right.
%\label{m^M}
%\end{equation}
It is well known (see, for example \cite{D}) that under the above restrictions on the function $\theta,$ the uncertainty constants of  for the Meyer scaling and wavelet function are bounded.

\texttt{The uncertainty constant} of $f$ is the functional $\Delta_{f}\Delta_{\widehat{f}}$ such that
$$
\begin{array}{ll}
\Delta_{f}^2:=\|f\|^{-2}_{L^2(\mathbb{R})}\int_{\mathbb{R}}(t-t_{0f})^2|f(t)|^2\,dt, &
\Delta_{\widehat{f}}^2:=\|\widehat{f}\|^{-2}_{L^2(\mathbb{R})}\int_{\mathbb{R}}(\w-\w_{0\widehat{f}})^2|\widehat{f}(\w)|^2\,d\w, \\
\phantom{1} & \phantom{1}\\
t_{0f}:=\|f\|^{-2}_{L^2(\mathbb{R})}\int_{\mathbb{R}}t|f(t)|^2\,dt, &
\w_{0\widehat{f}}:=\|\widehat{f}\|^{-2}_{L^2(\mathbb{R})}\int_{\mathbb{R}}
\w|\widehat{f}(\w)|^2\,d\w. \\
\end{array}
$$
The terms $\Delta_{f},$ $\Delta_{\widehat{f}},$ $t_{0f},$ and $\w_{0\widehat{f}}$ are
called \texttt{the time radius}, \texttt{the frequency radius}, \texttt{the time centre}, and \texttt{the frequency centre} of the function $f$ respectively.

The numbers $\pm e^{i\bar{\w}}$ are called \texttt{the pair of symmetric roots} of the mask $m$ if $m(\bar{\w})=m(\bar{\w}+\pi)=0$. A set $B:=\left\{b_1,\ldots,b_n\right\}$ of distinct complex numbers is called \texttt{cyclic} if $b_{j+1}=b_j^2,$ for $j=1,\ldots,n$ and $b_{n+1}=b_1.$ A cyclic set  $B$ is called \texttt{the cycle of the mask} $m$ if $m(\w+\pi)=0$ for all $\w$ such that $\exp(i\w)=b_j$ for some $j=1,\ldots,n.$ \texttt{The trivial cycle} is the set $\{1\}$. A mask is called \texttt{pure} if it has neither pairs of symmetric zeros nor cycles.
The following result gives a necessary and sufficient condition for integer shifts $\varphi(\cdot+k),$ $k \in \mathbb{Z}$ of a scaling function $\varphi$ to be stable (i.e., to form a Riesz basis).
\begin{prop}\cite[Corollary~3.4.15]{NPS}\label{Riesz} 
Integer shifts of a scaling function are stable (i.e., form a Riesz basis)
iff
corresponding mask has neither pairs of symmetric zeros nor nontrivial cycles.
\end{prop}

\texttt{The H\"{o}lder exponent} $\alpha_f$  of a function $f$ defined on some closed interval $[a,b]$ is 
$$
\alpha_f:=k+\sup_{\beta \in \mathbb{R}}\bigl\{\beta \in \mathbb{R}\bigl||f^{(k)}(x_1)-f^{(k)}(x_2)|\leq C_{\beta}|x_1-x_2|^{\beta}, \ x_1,x_2\in [a,b]\bigr.\bigr\}
,$$
 where $k:=\max_{h\in\mathbb{Z}}\bigl\{h\bigl|f\in C^{h}[a,b]\bigr.\bigr\}.$
Another characteristic of the smoothness of $f$ is
$$
\theta_{\widehat{f}}:=\sup_{\beta\in\mathbb{R}}\left\{\beta\in \mathbb{R}\left||\widehat{f}(\w)|\leq C(|\w|+1)^{-\beta}\right.\right\}.
$$
The smoothness characteristics we introduced are known to satisfy the inequality 
$\theta_{\widehat{f}}-1\leq \alpha_f\leq \theta_{\widehat{f}}.$
By $\theta(m)$ we mean $\theta_{\widehat{\varphi}},$ where $\varphi$ is the scaling function corresponding to the mask  $m$.   The following
result can be used for finding $\theta(m).$

\begin{prop}\cite[Lemma~7.4.2 and Proposition~7.4.4]{NPS}\label{smoo}
Suppose that some mask  $m$ is represented as
$m(\w)=\left(\cos \frac \w2\right)^{L+1}m_c(\w),$
where $m_c$ is a pure mask.
Then 
$\theta(m)=L+1+\theta(m_c)$
and
$\theta(m_c)=\lim_{k\rightarrow\infty}\theta_k,$
where
\begin{equation}
    \theta_k:=-\frac 1k \log_2\|m_c(\w)\cdots m_c(2^{k-1}\w)\|_{\infty}.
\label{thet}
\end{equation}
\end{prop}

\section{Basic construction and conditions for a linear method of summation}
Let us introduce a non-orthogonal mask of a new wavelet function. 
It is defined as the following $2\pi$-periodic trigonometric polynomial
\begin{equation}
\label{ml}
	m_l(\w):=\left(\cos\frac \w2\right)^{2l}
	 \frac{u_{n(l)}(m^M_l,\w)}{u_{n(l)}(m^M_l,0)},
\end{equation}
where
$$
m^M_l(\w):=\frac{m^M(\w)}{\left(\cos\frac \w2\right)^{2l}},\; l\in\mathbb{N},
$$
$m^M$  is a fixed Meyer mask, and trigonometric polynomial
$u_{n(l)}(m^M_l,\cdot)$ is defined by a fixed linear method of summation for the function 
$m^M_l.$

Since $m_l$ is a trigonometric polynomial and $m_l(0)=1,$  the infinite product $\prod_{j=1}^{\infty}m_l\left(\frac{\w}{2^j}\right)$ converges absolutely and uniformly on an arbitrary compact set. (If an infinite
product is equal to zero we assume that it converges.)   Thus the function $m_l$ is a mask for a stable,
but not orthogonal scaling function $\varphi_l,$ the Fourier transform of $\varphi_l$ is determined by the equality  
\begin{equation}
\label{ad1}
	\widehat{\varphi_l(\w)}=\prod_{j=1}^{\infty}m_l\left(\frac{\w}{2^j}\right).
\end{equation}
 The functions $\varphi_l(\cdot + k)$ for $k \in \mathbb{Z}$ form the Riesz basis in the closure of their linear span;
this claim is straight corollary of  Lemma \ref{pure} and Proposition \ref{Riesz}. 
%The function 
%$\varphi_l$ is called \textit{a non-orthogonal quasispline scaling  function}.
From  estimation (\ref{ocenka}) it follows that the orthogonalizing factor 
\begin{equation}
\label{ad2}
\Phi_l(\w):=\sum_{k\in \mathbb{Z}}\left|\widehat{\varphi_l}(\w+2\pi k)\right|^2
\end{equation}
is well defined.
Using the function $
\Phi_l$ we define the Fourier transform of the orthogonal scaling function 
\begin{equation}
\label{ad3}
\widehat{\varphi_l^{\bot}}(\w):=\widehat{\varphi_l}(\w)\Phi_l^{-0,5}(\w),
\end{equation}
the orthogonal mask 
\begin{equation}
\label{ad4}
m_l^{\bot}(\w):=m_{l}(\w)\Phi^{0,5}(\w)\Phi^{-0,5}(2\w),
\end{equation}
and, finally, the Fourier transform of the wavelet function 
\begin{equation}
\label{qsw}	\widehat{\psi^{\bot}_l}(\w):=e^{\frac{-i\w}{2}}\overline{m_l^{\bot}}\left(\frac{\w}{2}+\pi\right)
	\widehat{\varphi^{\bot}_l}\left(\frac{\w}{2}\right).
\end{equation}
%So a \textit{quasispline wavelet function} is defined in Fourier domain by (\ref{qsw}). A trigonometric polynomial defined by (\ref{ml}) is called \textit{a non-orthogonal quasispline mask}.

\begin{defi}
\label{qswdf}
By \texttt{a quasispline wavelet function} we mean the function 
$\psi^{\bot}_l,$ where the Fourier transform  $\widehat{\psi^{\bot}_l}$ is defined by (\ref{qsw}) and a non-orthogonal mask is defined by (\ref{ml}). The functions 
$\varphi_l^{\bot},$ 
$m_l^{\bot},$ $\varphi_l,$ $m_l$  defined by (\ref{ad3}),(\ref{ad4}),(\ref{ad1}), and (\ref{ml}) respectively are called
\texttt{a quasispline scaling function}, \texttt{a quasi\-spline mask}, \texttt{a non-orthogonal quasispline scaling function}, and \texttt{a non-orthogonal quasispline mask} respectively.  
%
%A wavelet function defined in the Fourier domain by (\ref{qsw}), where a corresponding non-orthgonal mask is defined by (\ref{ml}), is called a \textit{quasispline wavelet function}.
\end{defi}
So for any fixed Meyer mask and for any fixed linear method of summation we get the sequence $(\psi_l)_{l \in \mathbb{N}}$ of quasispline wavelet functions, and the symbol $l$ is a smoothness parameter (see Theorem \ref{Smo}).

In the remaining part of the article the following main Theorem will be proven.
\begin{teo}
\label{mainteo}
Suppose that $\psi^{\bot}_l$ ($\varphi^{\bot}_l$) is a quasispline wavelet (scaling) function.
Then
\begin{enumerate}
	\item the functions 
$\varphi_l^{\bot}$ and $\psi_l^{\bot}$ decay exponentially at infinity (Theorem \ref{expp});
	\item the functions $\widehat{\varphi_l^{\bot}}$ and $\widehat{\psi_l^{\bot}}$ decay as $O(\w^{-l})$ at infinity, namely the H\"{o}lder exponents $\alpha_{\varphi_l^{\bot}}$ and $\alpha_{\psi_l^{\bot}}$ of the functions satisfy the inequalities  
	$$
2l-1+\log_2\left(\frac {c}{1+\ve(l)} \right)\leq \alpha_{\varphi_l^{\bot}}\leq 2l,
\quad
2l-1+\log_2\left(\frac  {c}{1+\ve(l)} \right)\leq \alpha_{\psi_l^{\bot}}\leq 2l
$$
(Theorem \ref{Smo});
	\item the uncertainty constants 	$\Delta^2_{\varphi_l^{\bot}} \Delta^2_{\widehat{\varphi_l^{\bot}}}$  ($\Delta^2_{\psi_l^{\bot}} \Delta^2_{\widehat{\psi_l^{\bot}}}$) 
of the quasispline scaling (wavelet)  functions 
	 $\widehat{\varphi_l^{\bot}}$ ($\psi_l^{\bot}$) tend to those of the Meyer scaling (wavelet) function,
	namely	
	$$
|\Delta^2_{\widehat{\varphi_l^{\bot}}}-\Delta^2_{\widehat{\varphi^M}}|=
O\left(\max\{\mu(l),\, (4e^{2\w_0})^{-2l+4\log_2 \frac {1+\ve(l)}{c}}\}\right),
$$
$$
|\Delta^2_{\varphi_l^{\bot}}-\Delta^2_{\varphi^M}|=
O\left(\max\{\mu(l),\, l C_0^{-l+2\log_2 \frac {1+\ve(l)}{c}}\}\right),
$$	
$$
|\Delta^2_{\widehat{\psi_l^{\bot}}}-\Delta^2_{\widehat{\psi^M}}|=
O\left(\max\{\mu(l),\, l C_0^{-l+2\log_2 \frac {1+\ve(l)}{c}}\}\right),
$$
$$
|\Delta^2_{\psi_l^{\bot}}-\Delta^2_{\psi^M}|=
O\left(\max\{\mu(l),\, (4 e^{2\w_0})^{-l+2\log_2 \frac {1+\ve(l)}{c}}\}\right)
$$	
as $l \to \infty,$	
where
$C_0:=\frac{32 \pi^2 e^{2\w_0}}{27},$ others parameters 
are defined by (\ref{param}) (Theorems \ref{cfr}, \ref{cti}  and \ref{cfrtiw}).
\end{enumerate}
\end{teo}

For the quasispline wavelet function $\psi^{\bot}_l$ to satisfy  Theorem \ref{mainteo} and therefore to solve problem (\ref{prob}) it is sufficient to have three following conditions for the polynomials   $u_{n(l)}(m^M_l,\cdot).$ 

\begin{teo}\label{cond}
Suppose that there exists a sequence $n(l)$ for $l \in \mathbb{N}$ such that
\begin{equation}
\label{con1}
	\|u_{n(l)}(m^M_l,\cdot)-m^M_l\|_{C}=:\alpha(l)=o(l^{-1}) \mbox{ as } l \to \infty,
\end{equation}
\begin{equation}
\label{con2}
	\|u_{n(l)}((m^M_l)',\cdot)-(m^M_l)'\|_{C}=:\gamma(l)=o(1) \mbox{ as }  l \to \infty,
\end{equation}
\begin{equation}\label{con3}
	u_{n(l)}(m^M_l,\pi)\neq 0,
\end{equation}
then the corresponding quasispline scaling (\ref{ad3}) and wavelet  (\ref{qsw}) functions satisfy the conditions of Theorem \ref{mainteo} 
\end{teo}
De la Vallee Poussin means satisfy these conditions (for the proof see \cite{LSMZ} p. 460, p.465, and p. 461 respectively).
 
By definition, put
$u_{l}:=u_{n(l)}(m^M_l,\cdot)$ 
and $u_{1,l}:=u_{n(l)}((m^M_l)',\cdot),$
$u_{0,l}:=u_l/u_l(0).$

\section{Convergence of frequency radii for the scaling functions}
\begin{lem}
\label{Cml}
$
\|m_l-m^M\|_{C} \leq K \alpha(l)=o(l^{-1}) \mbox{ as } l \to \infty,
$
where $K:=\frac{\|u_l\|_{C}}{\inf_{k \geq l}|u_k(0)|}$ is bounded.
\end{lem}
\textbf{Proof.} 
Combining  (\ref{ml}) and (\ref{con1}) we get 
$$
\|m_l-m^M\|_{C}=\left\|\left(\cos (\w/2)\right)^{2l} \frac{u_l}{u_l(0)}-m^M\right\|_{C}
\leq \left\|\frac{u_l}{u_l(0)}-m^M_l\right\|_{C}
\leq 
$$ 
$$
\leq \left\|\frac{u_l}{u_l(0)}-u_l\right\|_{C}+
\left\|u_l-m^M_l\right\|_{C} \leq 
\frac{\|u_l\|_{C}}{\inf_{k\geq l}|u_k(0)|}\left\|u_l(0)-1\right\|_C+\alpha(l)
\leq
$$
$$
\leq 
\left(\frac{\|u_l\|_{C}}{\inf_{k \geq l}|u_k(0)|}+1\right)\alpha(l).
$$
Since $u_l(0)\to m^M_l(0)=1$ as $l \to \infty,$ it follows that
for some $l_0 \in \mathbb{N}$ $\inf_{k \geq l_0}|u_k(0)| \geq c >0,$ therefore
$\frac{\|u_l\|_{C}}{\inf_{k \geq l}|u_k(0)|}$ is bounded. $\Box$

From here we suppose that $l \geq l_0.$  
To simplify reading let us collect together notation of parameters using in estimations. So we get

\begin{equation}
\label{param}
	\begin{array}{l}
	\mu(l):=l \alpha(l)+\gamma(l),\quad 
	c:=\inf_{l \geq l_0}|u_l(0)|,\quad 
	\ve(l):=\frac{\alpha(l)}{\|m_l^M\|_C},\\
	\frac{\pi}{3}\leq \w_0 < \frac{\pi}{2} \ (\mbox{the parameter of the Meyer mask})
	\end{array}
\end{equation}

\begin{lem}
\label{Cm'l}
$
\|m'_l-(m^M)'\|_{C} =O(\mu(l)) \mbox{ as } l \to \infty.
$
\end{lem}

\textbf{Proof.}
Using Lemma \ref{Cml}, (\ref{u'}), and (\ref{con2}) we get
$$
\left|\left(\left(\cos\frac{\w}{2}\right)^{2l} u_l(\w)\right)'-(m^M)'(\w)\right|=
$$
$$
=\left|-l \left(\cos\frac{\w}{2}\right)^{2l-1} \sin \frac{\w}{2}\, u_l(\w)+\left(\cos\frac{\w}{2}\right)^{2l} u'_l(\w) -(m^M)'(\w)\right|=
$$ 
$$
=\left|-l \left(\cos\frac{\w}{2}\right)^{2l-1} \sin \frac{\w}{2} \, (m^M_l(\w)+u_l(\w)-m^M_l(\w))
+\right.
$$
$$
\left.
+\left(\cos\frac{\w}{2}\right)^{2l} \left((m^M_l)'(\w)+u_{1,l}(\w)-(m^M_l)'(\w)\right)-(m^M)'(\w)\right|=
$$
$$
=
\left|-l \tan \frac{\w}{2}\,  m^M(\w)-l \left(\cos\frac{\w}{2}\right)^{2l} \sin \frac{\w}{2} \, (u_l(\w)-m^M_l(\w))+
\right.
$$
$$
\left.+
\left(\cos\frac{\w}{2}\right)^{2l} \cdot \frac{(m^M)'(\w) \left(\cos\frac{\w}{2}\right)^{2l}+l \left(\cos\frac{\w}{2}\right)^{2l-1} \sin \frac{\w}{2}
\,  m^M(\w)}{\left(\cos\frac{\w}{2}\right)^{4l}}+
\right.
$$
$$
\left.+
\left(\cos\frac{\w}{2}\right)^{2l} \left(u_{1,l}(\w)-(m^M_l)'(\w)\right)
-(m^M)'(\w)\right|=
$$
$$
\left|-l \left(\cos\frac{\w}{2}\right)^{2l-1} \sin \frac{\w}{2}\, (u_l(\w)-m^M_l(\w))
+
\left(\cos\frac{\w}{2}\right)^{2l} \left(u_{1,l}(\w)-(m^M_l)'(\w)\right)
\right|=
$$
$$
=O(l \alpha(l)+\gamma(l)).
$$
For $m_l,$ we have
$$
\left|m'_l(\w)-(m^M)'(\w)\right|=
\left|\frac{\left(\left(\cos\frac{\w}{2}\right)^{2l} u_l(\w)\right)'}{u_l(0)}-(m^M)'(\w)\right|\leq
$$
$$
\leq
\left|\left(\left(\cos\frac{\w}{2}\right)^{2l} u_l(\w)\right)'\right|
\left|u^{-1}_l(0)-1\right|+
\left|\left(\left(\cos\frac{\w}{2}\right)^{2l} u_l(\w)\right)'-(m^M)'(\w)\right|=
$$
$$
=\left(\|(m^M)'(\w)\|_{C}
+O(l \alpha(l)+\gamma(l))\right)\frac{O(\alpha(l))}{c}+
O(l \alpha(l)+\gamma(l))=
$$
$$
=O(l \alpha(l)+\gamma(l)).\Box
$$

\begin{lem}\label{cphiC}
$
\|\widehat{\varphi_l}-\widehat{\varphi^M}\|_{C[a,\,b]}=O(\mu(l)) 
$
as
$l \to \infty$ for any $a<b,$ $a,b \in \mathbb{R}.$ Parameter $\mu(l)$
is defined by (\ref{param})
\end{lem}

\textbf{Proof.}
One can rewrite the proof of  the Lemma from \cite[Lemma 1]{LSMZ}. It is sufficient to change the notation $v_l$ by $u_l$ and so on and to use the conditions (\ref{con1}),
(\ref{con2}) instead of the property of the de la Vallee Poussin mean (see the formulas 
(4)-(7), (11), (12) \cite{LSMZ}) $\Box.$ 

\begin{lem}\label{cphiL}
$
\left\|\widehat{\varphi_l}-\widehat{\varphi^M}\right\|_{L^2(\mathbb{R})}=
O\left(\max\{\mu(l),\, (4e^{2\w_0})^{-l+2\log_2 \frac {1+\ve(l)}{c}}\}\right)$ 
as
$l\to\infty. 
$
Parameters 
are defined by (\ref{param}).
\end{lem}

 \textbf{Proof.}
We claim that there exists a function   $\xi$ such that
$\xi\in L^2(\mathbb{R})$  and
$\left|\widehat{\varphi_l}(\w)\right|\leq\xi(\w).$
The construction of the majorant can be
rewritten with a inessential changes of notation from \cite[Lemma 2]{LSMZ}.
So write the results.
Denote 
\begin{equation}\label{phi0}
\widehat{\varphi_{0,l}}(\w):=\prod_{j=1}^{\infty}
\frac{u_l\left(\w2^{-j}\right)}{u_l(0)}.
\end{equation}
Then under the assumption $|\w|\geq 1$ we have
\begin{equation}\label{phi0oc}
	\left|\widehat{\varphi_{0,l}}(\w)\right|\leq |\w|^{-2\theta(u_{0,l})}
e^{2\w_0(l+O(\mu(l)))}\leq |\w|^{2\log_2 \frac{1+\ve(l)}{c}}
e^{2\w_0(l+O(\mu(l)))}.
\end{equation}

So $|\widehat{\varphi_l}(\w)|$ are majorized by the functions
\begin{equation}
|\widehat{\varphi_l}(\w)|\leq\xi_l(\w):=
\left\{
\begin{array}{ll}
\left|\widehat{\varphi^M}(\w)\right|+O(\mu(l)), & |\w|\leq 4e^{2\w_0},\\
e^{O(\mu(l))}|\w|^{-l+2\log_2 \frac {1+\ve(l)}{c}}, & |\w|> 4e^{2\w_0}.
\end{array}
\right.
\label{ocenka}
\end{equation}

Thus the function  $\xi$ may be defined as
$$
\xi(\w):=\left\{
\begin{array}{ll}
\nu_1, & |\w|\leq 4e^{2\w_0},\\
\nu_2|\w|^{-l_1+2\log_2 \frac {1+\ve(l)}{c}}, & |\w|> 4e^{2\w_0},
\end{array}
\right.
$$
where $\nu_1$ and $\nu_2$ are  
constants, $\nu_1,\nu_2>0$,
$l_1:=\max\{l_0,\; 2\log_2\frac {1+\ve(l)}{c} +2\}.$
Then the convergence follows  from the Lebesgue's dominated convergence Theorem and Lemma \ref{cphiC}. 

Let us estimate the rate of the convergence.
If $|\w|\geq 4e^{2\w_0},$ then $\widehat{\varphi^M}(\w)=0,$
so
$$
\left\|\widehat{\varphi_l}-\widehat{\varphi^M}\right\|^2_{L^2(\mathbb{R})}=
\int_{\mathbb{R}}\left|\widehat{\varphi_l}(\w)-\widehat{\varphi^M}(\w)\right|^2\,d\w=
\int_{|\w|<4e^{2\w_0}}+ \int_{|\w|\geq 4e^{2\w_0}}\leq
$$
$$
\leq
8e^{2\w_0}\|\widehat{\varphi_l}-\widehat{\varphi^M}\|^2_{C[-4e^{2\w_0},\,4e^{2\w_0}]}+
e^{O(\mu(l))}
\int_{|\w|\geq 4e^{2\w_0}}|\w|^{-2l+4\log_2 \frac {1+\ve(l)}{c}}\,d\w=
$$
$$
=
8e^{2\w_0}\|\widehat{\varphi_l}-\widehat{\varphi^M}\|^2_{C[-4e^{2\w_0},\,4e^{2\w_0}]}+
\frac{e^{O(\mu(l))}(4e^{2\w_0})^{-2l+4\log_2 \frac {1+\ve(l)}{c}+1}}{2l-4\log_2 \frac {1+\ve(l)}{c}-1}.
$$
This completes the proof of Lemma \ref{cphiL}
$\Box.$
\begin{z}\label{cphiCR}
If we combine Lemma \ref{cphiC} and Lemma \ref{cphiL}, we get 
$\|\widehat{\varphi_l}-\widehat{\varphi^M}\|_{C(\mathbb{R})}=O\left(\max\{\mu(l),\, (4e^{2\w_0})^{-l+2\log_2 \frac {1+\ve(l)}{c}}\}\right).$
\end{z}

\begin{lem}\label{cPhi}
$
\left|\Phi_l(\w)-1\right|=O\left(\max\{\mu(l),\, (4e^{2\w_0})^{-2l+4\log_2 \frac {1+\ve(l)}{c}}\}\right)
$
as $l \to \infty.$ Parameters 
are defined by (\ref{param}).
\end{lem}

\textbf{Proof.}
Suppose $\w \in [-\pi,\,\pi].$
Since $\varphi^M$ is an orthogonal scaling function, we see that 
$\sum_{k \in \mathbb{Z}}|\widehat{\varphi^M}(\w +2 \pi \w)|^2=1.$
Taking into account (\ref{ocenka}), we define 
$k_0:=\left[2 e^{2 \w_0}/\pi+1/2\right].$
Hence
$$
\left|\Phi_l(\w)-1\right|=
\left|\sum_{k \in \mathbb{Z}}|\widehat{\varphi_l}(\w +2 \pi k)|^2-
\sum_{k \in \mathbb{Z}}|\widehat{\varphi^M}(\w +2 \pi k)|^2\right|
\leq
$$
$$
\leq
\sum_{k \in \mathbb{Z}}\left|\left(\widehat{\varphi_l}(\w +2 \pi k)\right)^2-
\left(\widehat{\varphi^M}(\w +2 \pi k)\right)^2\right|=
\sum_{|k|\leq k_0}+\sum_{|k|>k_0}.
$$
Using Lemma \ref{cphiC} we get
$$
\sum_{|k|\leq k_0}\leq 
(2k_0+1)(\sup_{|\w|\leq 4 e^{2 \w_0} }|\widehat{\varphi_l}(\w)-\widehat{\varphi^M}(\w)|+2\sup_{|\w|\leq 4 e^{2 \w_0} }|\widehat{\varphi^M}(\w)|)\times
$$
$$
\times
\sup_{|\w|\leq 4 e^{2 \w_0} }|\widehat{\varphi_l}(\w)-\widehat{\varphi^M}(\w)|\leq
O(\mu(l)).
$$
Since $\widehat{\varphi^M}=0$ as $|\w|\leq 4e^{2\w_0},$ (\ref{ocenka}), and the definition of $k_0,$ we obtain
$$
\sum_{|k|>k_0}\leq \sum_{|k|>k_0}
e^{O(\mu(l))}|\w+2\pi k|^{-2l+4\log_2 \frac {1+\ve(l)}{c}} =
O\left((4e^{2\w_0})^{-2l+4\log_2 \frac {1+\ve(l)}{c}}\right).
$$
Therefore,
$$
\left|\Phi_l(\w)-1\right|=O\left(\max\{\mu(l),\, (4e^{2\w_0})^{-2l+4\log_2 \frac {1+\ve(l)}{c}}\}\right) \Box.
$$

Now let us prove the convergence of the frequency radii for the scaling function.
 
\begin{teo}\label{cfr}
$
|\Delta^2_{\widehat{\varphi_l^{\bot}}}-\Delta^2_{\widehat{\varphi^M}}|=
O\left(\max\{\mu(l),\, (4e^{2\w_0})^{-2l+4\log_2 \frac {1+\ve(l)}{c}}\}\right)
$
as $l \to \infty.$ Parameters 
are defined by (\ref{param}).
\end{teo}

\textbf{Proof.}
Since the functions $\widehat{\varphi_l^{\bot}}$ and $\widehat{\varphi^M}$ are
even, then $\w_{0 \widehat{\varphi_l^{\bot}}}=\w_{0 \widehat{\varphi^M}}=0,$
where $\w_{0 \widehat{\varphi_l^{\bot}}},$ $\w_{0 \widehat{\varphi^M}}$ are the
frequency centers.

Taking into account Lemmas \ref{cphiC}, \ref{cPhi}, and the estimation (\ref{ocenka}) we have
$$
|\Delta^2_{\widehat{\varphi_l^{\bot}}}-\Delta^2_{\widehat{\varphi^M}}|=
\left|\int_{\mathbb{R}}\w^2 \left(\left(\widehat{\varphi_l^{\bot}}\right)^2(\w)-
\left(\widehat{\varphi^M}\right)^2(\w)\right)\, d\w \right|\leq
$$ 
$$
\leq \int_{\mathbb{R}} \w^2 \left|\frac{\left(\widehat{\varphi_l}\right)^2(\w)}{\Phi_l(\w)}-\left(\widehat{\varphi^M}\right)^2(\w)\right|\, d\w\leq \int_{|\w|<4e^{2 \w_0}}+\int_{|\w|\geq 4e^{2 \w_0}}\leq
$$
$$
\leq
16e^{4\w_0} \int_{|\w|<4e^{2 \w_0}}\left(\left(\widehat{\varphi_l}\right)^2(\w)\left|\frac{1}{\Phi_l(\w)}-1\right|+
\left|\left(\widehat{\varphi_l}\right)^2(\w)-\left(\widehat{\varphi^M}\right)^2(\w)\right|\right)\, d\w+
$$
$$
+
\int_{|\w|\geq 4e^{2 \w_0}}\w^2\left(\widehat{\varphi_l}\right)^2(\w)\frac{1}{\Phi_l(\w)}\,d\w
\leq
16e^{4\w_0} \left(\|\Phi_l-1\|_{C}\int_{|\w|<4e^{2 \w_0}}\frac{\left(\widehat{\varphi_l}\right)^2(\w)}{\Phi_l(\w)}\,d\w+\right.
$$
$$\left.
+\|\widehat{\varphi_l}-\widehat{\varphi^M}\|_{C[-4e^{2\w_0},\,4e^{2\w_0}]}\int_{|\w|<4e^{2 \w_0}}
|\widehat{\varphi_l}(\w)+\widehat{\varphi^M}(\w)|\,d\w\right)+
$$
$$
+\frac{2e^{O(\mu(l))} (4e^{2\w_0})^{-2l +4\log_2\frac{1+\ve(l)}{c}+3}}{\inf_{\w,\, l}\Phi_l(\w) \left(2l -4\log_2\frac{1+\ve(l)}{c}-3\right)}.
$$
From Lemmas \ref{cphiC} and \ref{cPhi} it follows that the integrals 
$$
\int_{|\w|<4e^{2 \w_0}}\frac{\left(\widehat{\varphi_l}\right)^2(\w)}{\Phi_l(\w)}\,d\w,
\quad 
\int_{|\w|<4e^{2 \w_0}}
|\widehat{\varphi_l}(\w)+\widehat{\varphi^M}(\w)|\,d\w
$$
are bounded. Hence
$$
|\Delta^2_{\widehat{\varphi_l^{\bot}}}-\Delta^2_{\widehat{\varphi^M}}|=
O\left(\max\{\mu(l),\, (4e^{2\w_0})^{-2l+4\log_2 \frac {1+\ve(l)}{c}}\}\right)+
O(\mu(l))+
$$
$$
+O((4e^{2\w_0})^{-2l+4\log_2 \frac {1+\ve(l)}{c}})=
O\left(\max\{\mu(l),\, (4e^{2\w_0})^{-2l+4\log_2 \frac {1+\ve(l)}{c}}\}\right) \Box.
$$

%\begin{z}\label{expphi}
%Since $m_l$  is a trigonometric polynomial, $\varphi_l$ is compactly supported. %Consequently, taking into account Lemma \ref{cPhi}, we deduce that $\varphi_l^{\bot}$ %decays exponentially at infinity.
%\end{z}

 \section{The growth of the smoothness and the exponential decaying.}
\begin{lem}
\label{pure}
The polynomial $u_{0,l}$ is a pure mask.
\end{lem}

\textbf{Proof.}
Let us use Proposition \ref{Riesz}.
Recall that $u_{0,l}=u_l/u_l(0)$.
By the condition (\ref{con1}) and the inequality $\pi/3 \leq \w_0 < \pi/2,$ where $\w_0$ is a parameter of the Meyer mask, we have 
$\sup_{[- \pi/3, \, \pi/3]}|u_l-m_l^M|
=\sup_{[- \pi/3, \, \pi/3]}
|u_l-(\cos \cdot /2)^{-2l}|=O(\alpha(l))$ as $l \to \infty.$
Hence
$u_l(\w)\neq 0$ on the interval $\w \in [- \pi/3, \, \pi/3].$
Therefore 
the polynomial $u_l$ has no pair of symmetric zeros.
If 
$B:=\left\{b_1,\ldots,b_n\right\}$ is a cyclic set and 
$b_1=r e^{i \xi},$ then $r=1,$ $\xi=\frac{2\pi k}{2^n-1}.$ 
If we suppose that 
$B$ is  a nontrivial cycle of the mask $u_l$ then the set $\pi+\frac{2\pi k}{2^n-1}$
has to be roots of  $u_l.$ But it does not hold true because of 
$u_l(\w)\neq 0$ on the interval $\w \in [- \pi/3, \, \pi/3].$
Finally, the condition $u_l(\pi)\neq 0$ is postulated in (\ref{con3}). Then  $u_l$ has no the trivial cycle
$\Box.$ 

Using Lemma \ref{pure} one can apply Proposition  \ref{smoo} to estimate smoothness of 
the non-orthogonal quasispline scaling function $\varphi_l.$

\begin{lem}\label{smophi}
The following inequality holds true
$2l-1+\log_2\left(\frac {c}{1+\ve(l)}\right)\leq \alpha_{\varphi_l}\leq 2l.$
Parameters 
are defined by (\ref{param}).
\end{lem}

\textbf{Proof.}
If we recall (\ref{con1}) and $c=\inf_{l\geq l_0}u_l(0),$ we get
$$
\sup_{\w}\left|u_{0,l}(\w)\right|\leq
\frac{\sup_{\w}|u_l(\w)|}{c}\leq
\frac{(1+\ve(l))\sup_{\w}|m^M_l(\w)|}{c}\leq 
\sup_{\w}|f_{0,l}(\w)|
$$
for $\ve(l):=\alpha(l)/\|m_l^M\|_C \to 0$ as $l \to \infty,$
where $f_{0,l}$ is even $2\pi$-periodic function and
$f_{0,l}(\w):=(1+\ve(l))(\cos \w/2)^{-2l}/c$ for $0\leq \w\leq \w_1$ and 
$f_{0,l}(\w):=0$ for $\w_1 < \w \leq \pi.$ 
So we get
$\theta_k(u_{0,l})\geq \theta_k(f_{0,l}).$

The definition of $f_{0,l}$ yields 
$$
\left\|f_{0,l}(\w)\cdots f_{0,l}(2^{k-1}\w)\right\|_{\infty}=f_{0,l}(\w_1)\cdots f_{0,l}(2^{-k+1}\w_1)=
$$
$$
=
\left(\cos\frac{\w_1}{2}\cdots\cos\frac{\w_1}{2^k}\right)^{-2l}\left(\frac {1+\ve(l)}{c} \right)^k.
$$
Then using Proposition \ref{smoo} we have
$$
\theta_k(f_{0,l})=-\frac 1k \log_2\left(\frac {1+\ve(l)}{c} \right)^k- 2l\log_2
\left|\cos\frac{\w_1}{2}\cdots\cos\frac{\w_1}{2^k}\right|^{-\frac 1k}
\to
\log_2\left(\frac {c}{1+\ve(l)} \right)
$$
as $k \to \infty.$
Passing to the limit, we use the identity
 $\prod_{j=1}^{\infty}\cos\frac{\w}{2^j}=\frac{\sin \w}{\w}.$
Therefore 
$\theta(u_{0,l})\geq \log_2\left(\frac {c}{1+\ve(l)} \right)$.  For $u_{0,l}$ the multiplicity of the trivial cycle is equal to $2l.$   Hence
 $2l-1+\log_2\left(\frac {c}{1+\ve(l)}\right)\leq \alpha_{\varphi_l}.$
By  definition of the norm $\|\cdot\|_{\infty}$ we have
$\|u_{0,l}(\w)\dots u_{0,l}(2^{k-1}\w)\|_{\infty}\geq u_{0,l}(0)\dots u_{0,l}(2^{k-1} \cdot 0)=1.$
Therefore Proposition \ref{smoo} yields
$\theta_k(u_{0,l})\leq 0,$
then
$\theta(u_{0,l})\leq 0,$
thus
$\alpha_{\varphi_l}\leq 2l.$
Finally, we obtain
$2l-1+\log_2\left(\frac {c}{1+\ve(l)} \right)\leq \alpha_{\varphi_l}\leq 2l. \Box$

Lemma \ref{cPhi} allows to extend the estimation of the smoothness to the orthogonal scaling and wavelet functions.
\begin{teo}\label{Smo}
The following inequalities hold true
$$
2l-1+\log_2\left(\frac {c}{1+\ve(l)} \right)\leq \alpha_{\varphi_l^{\bot}}\leq 2l,
\quad
2l-1+\log_2\left(\frac  {c}{1+\ve(l)} \right)\leq \alpha_{\psi_l^{\bot}}\leq 2l.
$$
Parameters 
are defined by (\ref{param}).
\end{teo}
 
 \textbf{Proof.}
It is sufficient to prove
$
\theta_{\widehat{\varphi_l}}=\theta_{\widehat{\varphi_l^{\bot}}}=\theta_{\widehat{\psi_l^{\bot}}}.
$
Using Lemma \ref{cPhi} we get $0<c_1\leq \Phi_l(\w)\leq c_2<\infty.$
Therefore $c_2^{-0.5}|\widehat{\varphi_l}|\leq |\widehat{\varphi_l^{\bot}}|\leq c_1^{-0.5}|\widehat{\varphi_l}|.$
Thus taking into account the definition of $\theta_{\widehat{f}}$ we get 
$
\theta_{\widehat{\varphi_l}}=\theta_{\widehat{\varphi_l^{\bot}}}.$

Then the application of (\ref{qsw}) yields
$$
|\widehat{\psi_{l}^{\bot}}(\w)|=\left|m_l\left(\frac{\w}{2}+\pi\right)
\Phi_l^{0,5}\left(\frac{\w}{2}+\pi\right)\Phi_l^{-0,5}\left(\w+2\pi\right)
\widehat{\varphi_l}\left(\frac{\w}{2}\right)\Phi_l^{-0,5}\left(\frac{\w}{2}\right)\right|.
$$
There exists an arbitrary large $\w$ (for example, $\w\in[-2\w_0+2\pi(2k-1), 2\w_0+2 \pi(2 k-1)],$  $k\in \mathbb{Z}$) such that 
$1-\alpha(l)\leq m_l(\w/2+\pi)\leq 1+\alpha(l).$ Therefore for given $\w$ we have 
$(1-\alpha(l)) c_1^{0.5} c_2^{-1}|\widehat{\varphi_{l}}(\w/2)|\leq|\widehat{\psi_{l}^{\bot}}(\w)|\leq (1+\alpha(l)) c_2^{0.5} c_1^{-1}|\widehat{\varphi_{l}}(\w/2)|.$
Finally, again taking into account the definition of $\theta_{\widehat{f}},$ we get
$
\theta_{\widehat{\varphi_l}}=\theta_{\widehat{\psi_l^{\bot}}}
$
$\Box.$

Lemma \ref{cPhi} also allows to deduce exponential decay of the orthogonal scaling function 
$\varphi_l^{\bot}$ and the wavelet fucntion $\psi_l^{\bot}.$ 
\begin{teo}\label{expp}
The functions 
$\varphi_l^{\bot}$ and $\psi_l^{\bot}$ decay exponentially at infinity.
\end{teo}

\textbf{Proof.}
Since $m_l$  is a trigonometric polynomial, $\varphi_l$ is compactly supported. Consequently taking into account Lemma \ref{cPhi}, we deduce that $\varphi_l^{\bot}$ decays exponentially at infinity. So $\varphi_l^{\bot}(t)=O(e^{-\beta_2|t|}),$ $\beta_ 2> 0.$
Fix $l\in \mathbb{N}.$
The application of equality (\ref{qsw}) yields
$\psi_l^{\bot}(t)=\sum_{k \in \mathbb{Z}}(-1)^k h_{-k+1}\varphi_l^{\bot}(2t-k),$
where $h_k$ are the Fourier coefficients of the function $m_l^{\bot}.$ As 
$m_l^{\bot}$ is a rational trigonometric function and the denominator does not equal to $0$ as $\w \in \mathbb{R},$ then 
$h_k=O(e^{-\beta_1|k|}),$ $\beta_1 > 0.$  
Therefore, we have
$$
|\psi^{\bot}_l(t)|=\left|\sum_{k\in \mathbb{Z}}(-1)^k h_{-k+1}\varphi^{\bot}_l(2 t-k)\right|\leq
\sum_{k\in \mathbb{Z}}\Bigl|h_{-k+1}\varphi^{\bot}_l(2 t-k)\Bigr|\leq
$$
$$
\leq A \sum_{k\in \mathbb{Z}}e^{-\beta_2|2t-k|-\beta_1|-k+1|}, 
$$
where $A$ is a constant.
The application of the property of module and geometric series yields
$$
\sum_{k\in \mathbb{Z}}e^{-\beta_2|2t-k|-\beta_1|-k+1|}=
$$
$$
=
\frac{e^{\pm\beta_2-2\beta_2|t|}}{1-e^{-\beta_1-\beta_2}}+
\frac{e^{\pm\beta_2}}{e^{\beta_2-\beta_1}-1}\Bigl(e^{\beta_2|-2t+[2t]|-\beta_1|[2t]|}-e^{-2\beta_2|t|}\Bigr)+
\frac{e^{\kappa+\beta_2|2t-[2t]|-\beta_1|[2t]|}}{1-e^{-\beta_1-\beta_2}},
$$
where $\kappa=-\beta_2$ as $t\geq 0$ and $\kappa=-\beta_1$ as $t<0.$
Therefore
$\psi^{\bot}_l=O\Bigl(e^{-\max\{\beta_1,\,\beta_2\}|2t|}\Bigr).$ $\Box$

\section{Convergence of time radii for the scaling functions}
\begin{lem}\label{cphi'C}
For any $-\infty<a<b<\infty$
it holds true
$
\left\|\widehat{\varphi_l}'-\widehat{\varphi^M}'\right\|_{C[a,\,b]}=
O\left(\mu(l)\right)$ 
as
$
l\to\infty.
$
Parameter $\mu(l)$ 
is defined by (\ref{param}).
\end{lem}

\textbf{Proof.}
Using the definition of 
$\widehat{\varphi_l}$ we get
$$
\left|\widehat{\varphi_l}'(\w)-\widehat{\varphi^M}'(\w)\right|=
\left|\left(\prod_{j=1}^{\infty}m_l\left(\frac{\w}{2^j}\right)\right)'-
\left(\prod_{j=1}^{\infty}m^M\left(\frac{\w}{2^j}\right)\right)'\right|
=$$
$$
=
\left|\sum_{j_0=1}^{\infty}2^{-j_0}\left(m_l'\left(\frac{\w}{2^{j_0}}\right)
\prod_{j=1, j\neq j_0}^{\infty}m_l\left(\frac{\w}{2^j}\right)-
(m^M)'\left(\frac{\w}{2^{j_0}}\right)
\prod_{j=1, j\neq j_0}^{\infty}m^M\left(\frac{\w}{2^j}\right)
\right)\right|\leq
$$
$$
\leq
\sum_{j_0=1}^{\infty}2^{-j_0}\left(
\left|m_l'\left(\frac{\w}{2^{j_0}}\right)-(m^M)'\left(\frac{\w}{2^{j_0}}\right)\right|
\prod_{j=1, j\neq j_0}^{\infty} \left|m^M\left(\frac{\w}{2^j}\right)\right|+\right.
$$
$$
\left.+
\left|m_l'\left(\frac{\w}{2^{j_0}}\right)\right|
\left|\prod_{j=1, j\neq j_0}^{\infty}m_l\left(\frac{\w}{2^j}\right)-\prod_{j=1, j\neq j_0}^{\infty}m^M\left(\frac{\w}{2^j}\right)\right|
\right).
$$
From Lemma \ref{Cm'l} it follows that
$
\left|m_l'\left(\frac{\w}{2^{j_0}}\right)-(m^M)'\left(\frac{\w}{2^{j_0}}\right)\right|
=O(\mu(l))
$
and
$
\left|m_l'\left(\frac{\w}{2^{j_0}}\right)\right|=\overline{M}+O(\mu(l)),
$
where 
$\overline{M}:=\left\|(m^M)'\right\|_C.$
Since $|m^M|\leq 1,$ we have
$
\left|\prod_{j=1, j\neq j_0}^{\infty}m^M\left(\frac{\w}{2^j}\right)\right|\leq 1.
$

Taking into account Lemma \ref{Cml} and the definition of 
$\widehat{\varphi_l},$ we obtain
$$
\left|\prod_{j=1, j\neq j_0}^{\infty}m_l\left(\frac{\w}{2^j}\right)-\prod_{j=1, j\neq j_0}^{\infty}m^M\left(\frac{\w}{2^j}\right)\right| \leq
\left|\prod_{j=1}^{j_0-1}m_l\left(\frac{\w}{2^j}\right)-\prod_{j=1}^{j_0-1}m^M\left(\frac{\w}{2^j}\right)\right| \left|\widehat{\varphi_l}\left(\frac{\w}{2^{j_0}}\right)\right|+
$$
$$
+
\left|\prod_{j=1}^{j_0-1}m^M\left(\frac{\w}{2^j}\right)\right|
\left|\widehat{\varphi_l}\left(\frac{\w}{2^{j_0}}\right)-\widehat{\varphi^M}\left(\frac{\w}{2^{j_0}}\right)\right|
$$
Using (\ref{con1}) and the property of the Meyer mask $m^M\leq 1$ we get
$$
\left|\prod^{j_0-1}_{j=1}m_l\left(\frac{\w}{2^j}\right)-
\prod^{j_0-1}_{j=1}m^M\left(\frac{\w}{2^j}\right)\right|\leq
$$
$$
\leq
\left|m_l\left(\frac{\w}{2}\right)-m^M\left(\frac{\w}{2}\right)\right|\prod^{j_0-1}_{j=2}m^M\left(\frac{\w}{2^j}\right)+
\left|m_l\left(\frac{\w}{2}\right)\right|\left|\prod^{j_0-1}_{j=2}m_l\left(\frac{\w}{2^j}\right)-
\prod^{j_0-1}_{j=2}m^M\left(\frac{\w}{2^j}\right)\right|\leq
$$
$$
\leq \|m_l-m^M\|_{C}+\left(1+\|m_l-m^M\|_{C}\right)\left|\prod^{j_0-1}_{j=2}m_l\left(\frac{\w}{2^j}\right)-
\prod^{j_0-1}_{j=2}m^M\left(\frac{\w}{2^j}\right)\right|.
$$
Reiterating the procedure $j_0-2$ times we obtain
$$
\left|\prod_{j=1}^{j_0-1}m_l\left(\frac{\w}{2^j}\right)-\prod_{j=1}^{j_0-1}m^M\left(\frac{\w}{2^j}\right)\right|=\left(1+O(\alpha(l))\right)^{j_0-1}-1.
$$
From Lemma \ref{cphiC} and the definition of the Meyer scaling function it follows that
$\left|\widehat{\varphi_l}\left(\frac{\w}{2^{j_0}}\right)\right|=1+O(\mu(l))$
and $\left|\widehat{\varphi_l}\left(\frac{\w}{2^{j_0}}\right)-\widehat{\varphi^M}\left(\frac{\w}{2^{j_0}}\right)\right|=O(\mu(l)).$
Finally, we note that $|m^M|\leq 1,$ therefore 
$\left|\prod_{j=1}^{j_0-1}m^M\left(\frac{\w}{2^j}\right)\right|\leq 1.$

Combining all the estimations together we obtain
$$
\left|\widehat{\varphi_l}'(\w)-\widehat{\varphi^M}'(\w)\right|\leq
O(\mu(l)) \sum_{j_0=1}^{\infty}2^{-j_0} +
\left(\overline{M}+O(\mu(l))\right) 
$$
$$
\left(O(\mu(l))\sum_{j_0=1}^{\infty}2^{-j_0}+\left(1+O(\mu(l))\right) \sum_{j_0=1}^{\infty}2^{-j_0}\left((1+O(\alpha(l)))^{j_0-1}-1\right)\right)=
$$
$$
=O(\mu(l))+\frac{O(\alpha(l))}{1-O(\alpha(l))}=O(\mu(l)).
$$
The next to last equality follows from the identity
$$
\sum_{j_0=1}^{\infty}2^{-j_0}\left((1+O(\alpha(l)))^{j_0-1}-1\right)=
$$
$$
=
\sum_{j_0=1}^{\infty}\left(\frac{1}{2}\left(\frac{1+O(\alpha(l))}{2}\right)^{j_0-1}-\frac{1}{2^{j_0}}\right)=\frac{O(\alpha(l))}{1-O(\alpha(l))} \Box.
$$

\begin{lem}\label{cphi'L}
$
\left\|\widehat{\varphi_l}'-\widehat{\varphi^M}'\right\|_{L_2(\mathbb{R})}=
O\left(\max\{\mu(l),\, l^{0.5} C_0^{-l+2\log_2 \frac {1+\ve(l)}{c}}\}\right)$
as
$l\to\infty,$
where $C_0:=\frac{32 \pi^2 e^{2\w_0}}{27}.$ Other parameters 
are defined by (\ref{param}).
\end{lem}

\textbf{Proof.}
We prove the Lemma in a similar manner as Lemma \ref{cphiL}. 
Let us find a majorant $\xi_1 \in L_2(\mathbb{R})$
for the function  
$\widehat{\varphi_l}'.$
From the definition of $\widehat{\varphi_l},$ (\ref{phi0}), and the identity 
$\sum_{j_0=1}^{\infty}2^{-j_0}  \tan\frac{\w}{2^{j_0+1}}=
\frac{2}{\w}-\cot\frac{\w}{2}
$
 it follows that
$$
\left(\widehat{\varphi_l}\right)'(\w)=
\left(\prod_{j=1}^{\infty}m_l\left(\frac{\w}{2^j}\right)\right)'=
\sum_{j_0=1}^{\infty}2^{-j_0}m_l'\left(\frac{\w}{2^{j_0}}\right)
\prod_{j=1, j\neq j_0}^{\infty}m_l\left(\frac{\w}{2^j}\right)=
$$
$$
=
\sum_{j_0=1}^{\infty}2^{-j_0}\left(l \left(\cos\frac{\w}{2^{j_0+1}}\right)^{2l-1}
\left(-\sin\frac{\w}{2^{j_0+1}}\right)\frac{u_l\left(\frac{\w}{2^{j_0}}\right)}{u_l(0)}
\prod_{j=1, j\neq j_0}^{\infty}\left(\cos\frac{\w}{2^{j+1}}\right)^{2l}\frac{u_l\left(\frac{\w}{2^{j}}\right)}{u_l(0)}+\right.
$$
$$
\left.+
\left(\cos\frac{\w}{2^{j_0+1}}\right)^{2l}
\frac{u_{1,l}\left(\frac{\w}{2^{j_0}}\right)}{u_l(0)}
\prod_{j=1, j\neq j_0}^{\infty}\left(\cos\frac{\w}{2^{j+1}}\right)^{2l}\frac{u_l\left(\frac{\w}{2^{j}}\right)}{u_l(0)}
\right)=
$$
$$
=\sum_{j_0=1}^{\infty}2^{-j_0} l \left(-\tan\frac{\w}{2^{j_0+1}}\right)
\prod_{j=1}^{\infty}\left(\cos\frac{\w}{2^{j+1}}\right)^{2l}
\prod_{j=1}^{\infty}\frac{u_l\left(\frac{\w}{2^{j}}\right)}{u_l(0)}+
$$
$$
+
\sum_{j_0=1}^{\infty}2^{-j_0}
\frac{u_{1,l}\left(\frac{\w}{2^{j_0}}\right)}{u_l(0)}
\prod_{j=1}^{\infty}\left(\cos\frac{\w}{2^{j+1}}\right)^{2l}
\prod_{j=1, j\neq j_0}^{\infty}\frac{u_l\left(\frac{\w}{2^{j}}\right)}{u_l(0)}=
$$
$$
=l \left(\cot\frac{\w}{2}-\frac{2}{\w}\right)\left(\frac{\sin \w/2}{\w/2}\right)^{2l}
\widehat{\varphi_{l,0}}(\w)+
$$
$$
+
\left(\frac{\sin \w/2}{\w/2}\right)^{2l}
\sum_{j_0=1}^{\infty}2^{-j_0}
\frac{u_{1,l}\left(\frac{\w}{2^{j_0}}\right)}{u_l(0)}
\prod_{j=1, j\neq j_0}^{\infty}\frac{u_l\left(\frac{\w}{2^{j}}\right)}{u_l(0)}
=:I_{1,l}(\w)+\left(\frac{\sin \w/2}{\w/2}\right)^{2l} I_{2,l}(\w).
$$
If $|\w| > 4 e^{2\w_0},$ then applying (\ref{phi0oc}) and (\ref{ocenka}) for the first item we have
$$
|I_{1,l}(\w)|=
\left|l \left(\cos\frac{\w}{2}-\frac{2\sin \w/2}{\w}\right)\left(\frac{2}{\w}\right)^{2l} \left(\sin \w/2\right)^{2l-1}
\widehat{\varphi_{l,0}}(\w)\right|\leq
$$
$$
\leq
C l e^{O(\mu(l))}  |\w|^{-l+2\log_2\frac{1+\ve(l)}{2}}.
$$
Let us estimate the second item. 
$$
I_{2,l}(\w)=
\sum_{j_0=1}^{\infty}2^{-j_0}
\frac{u_{1,l}\left(\frac{\w}{2^{j_0}}\right)}{u_l(0)}
\prod_{j=1}^{ j_0-1}\frac{u_l\left(\frac{\w}{2^{j}}\right)}{u_l(0)}
\widehat{\varphi_{l,0}}\left(\frac{\w}{2^{j_0}}\right).
$$
Using (\ref{phi0oc}) for  $|\w|\geq 4e^{2\w_0}$ we get
$\left|\widehat{\varphi_{l,0}}\left(\frac{\w}{2^{j_0}}\right)\right|
\leq
|\w 2^{-j_0}|^{-2 \theta(u_{0,l})} e^{2\w_0(l+O(\mu(l)))}.
$
Using the condition (\ref{con2}), the definition of the function $m^M_l,$ and 
the inequality $\pi/2 < \w_1 \leq 2\pi/3$ we get
$$
\left|u_{1,l}\left(\frac{\w}{2^{j_0}}\right)\right|\leq 
\left|(m^M_l)'\left(\frac{\w}{2^{j_0}}\right)\right|
+O(\gamma(l))\leq
$$
$$
\leq 
(m^M)'\left(\frac{\w}{2^{j_0}}\right)\left(\cos\frac{\w_1}{2^{j_0+1}}\right)^{-2l}+
l m^M\left(\frac{\w}{2^{j_0+1}}\right) \sin \frac{\w_1}{2^{j_0+1}}
\left(\cos\frac{\w_1}{2^{j_0+1}}\right)^{-2l-1}+
$$
$$
+O(\gamma(l))\leq
\left(\cos\frac{\w_1}{2^{j_0+1}}\right)^{-2l}\left(\overline{M}+l \tan \frac{\w_1}{2^{j_0+1}}\right)+O(\gamma(l))\leq
$$
$$
\leq (4/3)^l \left(\overline{M}+l \sqrt{3}\right)+O(\gamma(l)).
$$
Then
taking into account condition (\ref{con1}),
 the properties of the Meyer mask $|m^M|\leq 1,$ $m^M(\w)=0$ as $\w_1\leq|\w|\leq\pi,$
and the inequality $\pi/2 < \w_1 \leq 2 \pi/3$
we have
$$
\left|\prod_{j=1}^{ j_0-1}\frac{u_l\left(\frac{\w}{2^{j}}\right)}{u_l(0)}\right|\leq
\prod_{j=1}^{ j_0-1}\left(\frac{m^M(\w 2^{-j})}{\left(\cos \w 2^{-j-1}\right)^{2l}}
+\alpha(l)\right)
\leq
\prod_{j=1}^{ j_0-1}\left(\frac{1}{\left(\cos \w_1 2^{-j-1}\right)^{2l}}
+\alpha(l)\right)\leq
$$
$$
\leq
\prod_{j=1}^{ j_0-1}\frac{a}{\left(\cos \w_1 2^{-j-1}\right)^{2l}}=
a^{j_0-1}
\prod_{j=1}^{\infty}\left(\cos \w_1 2^{-j-1}\right)^{-2l}=
$$
$$
=a^{j_0-1}\left(\frac{\w_1/2}{\sin \w_1/2}\right)^{2l}
\leq a^{j_0-1} \left(\frac{\w_1}{\sqrt{2}}\right)^{2l}\leq
 a^{j_0-1} \left(\frac{2\pi}{3\sqrt{2}}\right)^{2l},
$$
where
$a$ is a majorant of the expression $1+\alpha(l)\left(\cos \w_1 2^{-j-1}\right)^{2l},$
so it can be chosen $a<1.5.$
 
Collecting the estimations we obtain for $I_{2,l}(\w)$
$$
I_{2,l}(\w)\leq
$$
$$
\leq
\sum_{j_0=1}^{\infty}2^{-j_0}
\frac{\left(\frac{4}{3}\right)^l \left(\overline{M}+l \sqrt{3}\right)+O(\gamma(l))}{1-\alpha(l)}
a^{j_0-1} \left(\frac{2\pi}{3\sqrt{2}}\right)^{2l}\left|\frac{\w}{ 2^{j_0}}\right|^{-2 \theta(u_{0,l})} e^{2\w_0(l+O(\mu(l)))}.
$$
Since
$\log_2\left(\frac {c}{1+\ve(l)} \right)\leq \theta(u_{0,l})\leq 0$ and
$a<1.5,$ we get $|\w|^{-2\theta(u_{0,l})}\leq |\w|^{2 \log_2\frac{1+\ve(l)}{c}}$ 
as $|\w|\geq 1,$
$2^{j_0 \theta(u_{0,l})}\leq 1,$ and $\sum_{j_0=1}^{\infty}2^{-j_0}a^{j_0-1}=(2-a)^{-1}.$
So
$$
I_{2,l}(\w)
\leq
\frac{ e^{O(\mu(l))}\left(\overline{M}+l \sqrt{3}+(3/4)^{l}O(\gamma(l))\right)}{(1-\alpha(l))(2-a)}
\left(\frac{8 \pi^2 e^{2 \w_0}}{27}\right)^l|\w|^{2 \log_2\frac{1+\ve(l)}{c}}
$$
Thus  we have for $|\w|>\frac{32 \pi^2 e^{2\w_0}}{27}$
$$
\left(\frac{\sin \w/2}{\w/2}\right)^{2l}
I_{2,l}(\w)\leq
\left(\sin \w/2\right)^{2l}
\frac{ e^{O(\mu(l))}\left(\overline{M}+l \sqrt{3}+(3/4)^{l}O(\gamma(l))\right)}{(1-\alpha(l))(2-a)}\times
$$
$$
\times
\left(\frac{32 \pi^2 e^{2 \w_0}}{27}\right)^l|\w|^{-2l+2 \log_2\frac{1+\ve(l)}{c}}
\leq C(l,\,\w_0) l |\w|^{-l+2 \log_2\frac{1+\ve(l)}{c}}, 
%= O\left(\left(\frac{32 \pi^2 e^{2\w_0}}{27}\right)^{-l+2 \log_2\frac{1+\ve(l)}{c}}\right),
$$
where
$C(l,\,\w_0):=
 e^{O(\mu(l))}\left(\overline{M}/l+ \sqrt{3}+l^{-1}(3/4)^{l}O(\gamma(l))\right)(1-\alpha(l))^{-1}(2-a)^{-1}$
is bounded with respect to the parameters $l$ and $\w_0.$ Put $C(l,\,\w_0)\leq A,$
$A$ is a constant.

So if  $|\w|>C_0:=\frac{32 \pi^2 e^{2\w_0}}{27},$ we can estimate  
$\left(\widehat{\varphi_l}\right)'$ as follows
$
\left|\left(\widehat{\varphi_l}\right)'(\w)\right|\leq A \, l\, 
|\w|^{-l+2 \log_2\frac{1+\ve(l)}{c}}.
$

Finally, using Lemma \ref{cphi'C} one can define the functions $\xi_{1,l}$ such that 

\begin{equation}\label{ocenka'}
	\left|\left(\widehat{\varphi_l}\right)'(\w)\right|\leq
	\xi_{1,l}(\w):=\left\{
	\begin{array}{ll}
	\left(\widehat{\varphi^M}\right)'(\w)+O(\mu(l)), & |\w|\leq \frac{32 \pi^2 e^{2\w_0}}{27},\\
	A \, l\, 
	|\w|^{-l+2 \log_2\frac{1+\ve(l)}{c}}, & |\w| \geq \frac{32 \pi^2 e^{2\w_0}}{27}.
	\end{array}
	\right.
\end{equation}

So the majorant $\xi_1$ is defined in the following way
$$
\xi_{1}(\w):=\left\{
\begin{array}{ll}
\nu'_1, & |\w|\leq \frac{32 \pi^2 e^{2\w_0}}{27},\\
\nu'_2 \, l\, 
|\w|^{-l_1+2 \log_2\frac{1+\ve(l)}{c}}, & |\w|\geq \frac{32 \pi^2 e^{2\w_0}}{27}.
\end{array}
\right.
$$

where $\nu'_1$ and $\nu'_2$ are  constants, $\nu'_1,\nu'_2>0$,
$l_1=\max\{l_0,\; 2\log_2\frac {1+\ve(l)}{c} +2\}$ is defined in the proof of Lemma \ref{cphiL}. 
Then the convergence follows  from the Lebesgue's dominated convergence Theorem and Lemma \ref{cphi'C}. 

Let us estimate the rate of the convergence.
If $|\w|\geq C_0,$ then $\widehat{\varphi^M}(\w)=0,$
so
$$
\left\|\widehat{\varphi_l}'-\widehat{\varphi^M}'\right\|^2_{L^2(\mathbb{R})}=
\int_{\mathbb{R}}\left|\widehat{\varphi_l}'(\w)-\widehat{\varphi^M}'(\w)\right|^2\,d\w=
\int_{|\w|<C_0}+ \int_{|\w|\geq C_0}\leq
$$
$$
\leq
2C_0\|\widehat{\varphi_l}'-\widehat{\varphi^M}'\|^2_{C[-C_0,\,C_0]}+
A^2 l^2
\int_{|\w|\geq C_0}|\w|^{-2l+4\log_2 \frac {1+\ve(l)}{c}}\,d\w=
$$
$$
=
2C_0\|\widehat{\varphi_l}'-\widehat{\varphi^M}'\|^2_{C[-C_0,\,C_0]}+
\frac{A^2 l^2(C_0)^{-2l+4\log_2 \frac {1+\ve(l)}{c}+1}}{2l-4\log_2 \frac {1+\ve(l)}{c}-1}=
$$
$$
=O\left(\max\{\mu^2(l),\, l C_0^{-2l+4\log_2 \frac {1+\ve(l)}{c}}\}\right).
$$
This completes the proof of Lemma \ref{cphi'L}
$\Box.$
\begin{z}\label{cphi'CR}
Using Lemma \ref{cphi'C} and Lemma \ref{cphi'L} we get
$$
\left\|\widehat{\varphi_l}'-\widehat{\varphi^M}'\right\|_{C(\mathbb{R})}
=O\left(\max\{\mu(l),\,l^{0.5} C_0^{-l+2\log_2 \frac {1+\ve(l)}{c}}\}\right).
$$
\end{z}

\begin{lem}\label{cPhi'}
$\|\Phi'_l\|_{C}= O\left(\max\{\mu(l),\,l^{0.5} (4 C_0 e^{2\w_0})^{-l+2\log_2 \frac {1+\ve(l)}{c}}\}\right)$ as $l \to \infty.$ Parameters 
are defined by (\ref{param}).
\end{lem}

\textbf{Proof.}
Taking into account the Definition $\Phi_l$ and the estimation (\ref{ocenka}) one can 
termwise differentiate the series, so
$$
\Phi'_l(\w)=\left(\sum_{k \in \mathbb{Z}}\left|\widehat{\varphi_l}(\w+2\pi k)\right|^2\right)'=\sum_{k \in \mathbb{Z}}\left(\left|\widehat{\varphi_l}(\w+2\pi k)\right|^2\right)'.
$$
Since the Meyer scaling function is compactly supported and satisfies the property
$\sum_{k \in \mathbb{Z}}\left|\widehat{\varphi^M}(\w+2\pi k)\right|^2\equiv 1,$
we get
$$
 \sum_{k \in \mathbb{Z}}\left(\left|\widehat{\varphi^M}(\w+2\pi k)\right|^2\right)'=
 \left(\sum_{k \in \mathbb{Z}}\left|\widehat{\varphi^M}(\w+2\pi k)\right|^2\right)'=
 (1)'=0.
$$
Suppose $|\w|\leq \pi,$ then we obtain
$$
|\Phi'_l(\w)| \leq 2\sum_{k \in \mathbb{Z}}\left|\widehat{\varphi_l}(\w+2\pi k)\widehat{\varphi_l}'(\w+2\pi k)-\widehat{\varphi^M}(\w+2\pi k)\widehat{\varphi^M}'(\w+2\pi k)\right|
\leq
$$
$$
\leq 2
\sum_{k \in \mathbb{Z}}\left|\widehat{\varphi_l}(\w+2\pi k)\right| \left|\widehat{\varphi_l}'(\w+2\pi k)-\widehat{\varphi^M}'(\w+2\pi k)\right|+
$$
$$
+
2 \sum_{k \in \mathbb{Z}}
\left|\widehat{\varphi^M}'(\w+2\pi k)\right|\left|\widehat{\varphi_l}(\w+2\pi k)-\widehat{\varphi^M}(\w+2\pi k)\right|=:2I_{3,l}(\w)+2I_{4,l}(\w).
$$
Using the parameter $k_0=\left[2 e^{2 \w_0}/\pi+1/2\right]$ 
defined in the proof of Lemma \ref{cPhi} we get
$$
I_{3,l}(\w)=\sum_{|k|\leq k_0}+\sum_{|k|> k_0}.
$$
Taking into account Lemma \ref{cphi'C}, for the first sum we have
$$
\sum_{|k|\leq k_0}\leq \|\widehat{\varphi_l}'-\widehat{\varphi^M}'\|_{C[-e^{2\w_0},\, e^{2\w_0}]}
\sum_{|k|\leq k_0}\left|\widehat{\varphi_l}(\w+2\pi k)\right|=O(\mu(l)).
$$
If we combine Remark \ref{cphi'CR} and the estimation (\ref{ocenka}), the second sum is
$$
\sum_{|k|> k_0} \leq
O\left(\max\{\mu(l),\, l^{0.5} C_0^{-l+2\log_2 \frac {1+\ve(l)}{c}}\}\right)
\sum_{|k|> k_0} |\w+2 \pi k|^{-l+2\log_2 \frac {1+\ve(l)}{c}}=
$$
$$
=O\left(\max\{\mu(l),\,l^{0.5} C_0^{-l+2\log_2 \frac {1+\ve(l)}{c}}\}\right)
(4 e^{2 \w_0})^{-l+2\log_2 \frac {1+\ve(l)}{c}}.
$$
%Thus
%$$
%I_{3,l}(\w)=O\left(\max\{\mu(l),\, C_0^{-l+2\log_2 \frac {1+\ve(l)}{c}}\}\right).
%$$
Estimate $I_{4,l}(\w).$
Since $\mbox{supp}\, \widehat{\varphi^M}=[-2\w_1,\,-2\w_0]\cup [2\w_0,\,2\w_1],$ then 
$\widehat{\varphi^M}(\w+2 \pi k)=0$ as $k>1.$
So for the sum $I_{4,l}(\w)$ we have
$$
I_{4,l}(\w)=\sum_{|k|\leq 1}
\left|\widehat{\varphi^M}'(\w+2\pi k)\right|\left|\widehat{\varphi_l}(\w+2\pi k)-\widehat{\varphi^M}(\w+2\pi k)\right|.
$$
Thus the application of Lemma \ref{cphiC} yields
$
I_{4,l}(\w)=O(\mu(l)).
$
Finally, for $\Phi'$ we get
$$
|\Phi'(\w)|\leq 2(I_{3,l}(\w)+I_{4,l}(\w))=O\left(\max\{\mu(l),\,l^{0.5} (4 C_0 e^{2\w_0})^{-l+2\log_2 \frac {1+\ve(l)}{c}}\}\right) \Box.
$$

Now let us prove the convergence of the time radii for the scaling function.

\begin{teo}\label{cti}
$
|\Delta^2_{\varphi_l^{\bot}}-\Delta^2_{\varphi^M}|=
O\left(\max\{\mu(l),\, l C_0^{-l+2\log_2 \frac {1+\ve(l)}{c}}\}\right)$
as $l \to \infty.$
Parameters 
are defined by (\ref{param}).
\end{teo}

\textbf{Proof.}
If the function $\widehat{\varphi}$ is real-valued, then
$\overline{\varphi}(t)=\varphi(-t).$ Hence
the function $|\varphi|^2$ is even. So
the time centre $t_{0 \varphi}=0.$
Then the square of the time radius
 $\Delta^2_{\varphi}=\int_{\mathbb{R}} t^2 |\varphi(t)|^2\,dt.$
 Using the property of the Fourier transform
 $\widehat{\varphi'}(\w)= i \w \widehat{\varphi}(\w)$ we obtain
 $\Delta^2_{\varphi}=(2\pi)^{-1}\int_{\mathbb{R}}|(\widehat{\varphi})'(\w)|^2\,d\w.$
 
Since the functions $\widehat{\varphi_l^{\bot}},$ $\widehat{\varphi^M}$ are
real-valued, then we have 
$t_{0 \varphi_l^{\bot}}=t_{0 \varphi^M}=0.$
So for the squares of the time radii we get 
$$
\Delta^2_{\varphi_l^{\bot}}=\frac{1}{2\pi}\int_{\mathbb{R}}|(\widehat{\varphi_l^{\bot}})'(\w)|^2\,d\w
\mbox{ and }
\Delta^2_{\varphi^M}=\frac{1}{2\pi}\int_{\mathbb{R}}|(\widehat{\varphi^M})'(\w)|^2\,d\w.
$$
Then we have
$$
\left|\Delta^2_{\varphi_l^{\bot}}-\Delta^2_{\varphi^M}\right|\leq
\frac{1}{2\pi}\int_{\mathbb{R}}\left|\left((\widehat{\varphi_l^{\bot}})'(\w)\right)^2-
\left((\widehat{\varphi^M})'(\w)\right)^2
\right|\,d\w\leq
$$
$$
\leq \frac{1}{2\pi}
\sup_{\w \in \mathbb{R}}\left|(\widehat{\varphi_l^{\bot}})'(\w)+
(\widehat{\varphi^M})'(\w)\right|
\int_{\mathbb{R}}\left|(\widehat{\varphi_l^{\bot}})'(\w)-
(\widehat{\varphi^M})'(\w)
\right|\,d\w
$$
Applying Lemmas \ref{cPhi}, \ref{cPhi'} and Remarks \ref{cphiCR}, \ref{cphi'CR} we establish the boundedness of the supremum factor
$$
\sup_{\w \in \mathbb{R}}\left|\widehat{\varphi_l^{\bot}}'(\w)+
\widehat{\varphi^M}'(\w)\right|\leq
\left\|\frac{\widehat{\varphi_l}'}{\Phi_l}\right\|_{C(\mathbb{R})}+
\left\|\frac{\Phi'_l\widehat{\varphi_l}}{\Phi^2_l}\right\|_{C(\mathbb{R})}+
\left\|\widehat{\varphi^M}'\right\|_{C(\mathbb{R})}\leq
$$
$$
\frac{\left\|\widehat{\varphi^M}'\right\|_{C(\mathbb{R})}+
\left\|\widehat{\varphi^M}'-\widehat{\varphi_l}'\right\|_{C(\mathbb{R})}}
{1-\left\|\Phi_l-1\right\|_{C}}+
$$
$$
+
\frac{\left\|\Phi'_l\right\|_{C} \left(\left\|\widehat{\varphi^M}\right\|_{C(\mathbb{R})}+
\left\|\widehat{\varphi^M}-\widehat{\varphi_l}\right\|_{C(\mathbb{R})}\right)}
{\left(1-\left\|\Phi_l-1\right\|_{C}\right)^2}+
\left\|\widehat{\varphi^M}'\right\|_{C(\mathbb{R})}=O\left(\left\|\widehat{\varphi^M}'\right\|_{C(\mathbb{R})}\right).
$$
Applying the same  Lemmas \ref{cPhi}, \ref{cPhi'} and Remarks \ref{cphiCR}, \ref{cphi'CR} we get the convergence to $0$ of the integral
$$
\int_{\mathbb{R}}\left|\widehat{\varphi_l^{\bot}}'(\w)-
\widehat{\varphi^M}'(\w)
\right|\,d\w \leq
\int_{\mathbb{R}}\left|\frac{\widehat{\varphi_l}'(\w)-
\widehat{\varphi^M}'(\w)\Phi_l(\w)}{\Phi_l(\w)}\right|\,d\w+
\int_{\mathbb{R}}\left|
\frac{\widehat{\varphi_l}(\w)\Phi'_l(\w)}{\Phi^2_l(\w)}\right|\,d\w\leq
$$
$$
\leq
\frac{1}{1-\left\|\Phi_l-1\right\|_{C}}
\int_{\mathbb{R}}\left|\widehat{\varphi_l}'(\w)-\widehat{\varphi^M}'(\w)\right|+
\left|\widehat{\varphi^M}'(\w)\right| \left|1-\Phi_l(\w)\right|\,d\w+
$$
$$
+\frac{\left\|\Phi'\right\|_{C}}{\left(1-\left\|\Phi_l-1\right\|_{C}\right)^2}
\int_{\mathbb{R}}\left|\widehat{\varphi_l}(\w)\right|\, d\w \leq
$$
$$
\leq
\frac{1}{1-\left\|\Phi_l-1\right\|_{C}}\left(
\int_{|\w|\leq C_0}\left|\widehat{\varphi_l}'(\w)-\widehat{\varphi^M}'(\w)\right| \,d\w+
\int_{|\w| > C_0}\left|\widehat{\varphi_l}'(\w)-\widehat{\varphi^M}'(\w)\right| \,d\w
\right)+
$$
$$
+\frac{\left\|1-\Phi\right\|_{C}}{1-\left\|\Phi_l-1\right\|_{C}}
\int_{\mathbb{R}}\left|\widehat{\varphi^M}'(\w)\right|\, d\w+
O\left(\left\|\Phi'\right\|_{C}\right)\leq
$$
$$
\leq
\frac{\left\|\widehat{\varphi_l}'-\widehat{\varphi^M}'\right\|_{C[-C_0,\,C_0]}}{1-\left\|\Phi_l-1\right\|_{C}}+\frac{2 A l C_0^{-l+2\log_2 \frac {1+\ve(l)}{c}+1}}{\left(1-\left\|\Phi_l-1\right\|_{C}\right)\left(l-2\log_2 \frac {1+\ve(l)}{c}-1\right)}+
$$
$$
+O\left(\left\|\Phi_l-1\right\|_{C}\right)+
O\left(\left\|\Phi'\right\|_{C}\right)=
O\left(\max\{\mu(l),\, l C_0^{-l+2\log_2 \frac {1+\ve(l)}{c}}\}\right)
\Box.
$$

In Lemmas \ref{cphi'C}-\ref{cPhi'} and Theorem \ref{cti} we apply without proof the formula
$\left(\prod_{j=1}^{\infty}m_l(\w 2^{-j})\right)'=\sum_{j_0=1}^{\infty}2^{-j_0}m'_l(\w 2^{-j_0})
\prod_{j=1, j\neq j_0}^{\infty} m_l(\w 2^{-j}).$ Establish it in the following

\begin{lem}\label{pdbp}
For any $a,\,b$ such that $-\infty <a<b<\infty$ we have
$$ 
\left\|\left(\prod_{j=1}^{\infty}m_l(\w 2^{-j})\right)'-\left(\prod_{j=1}^{n}m_l(\w 2^{-j})\right)'\right\|_{C[a,\,b]} \longrightarrow 0
$$
as $n \to \infty,$
where $\left(\prod_{j=1}^{\infty}m_l(\w 2^{-j})\right)'$ is the notation for the series

\noindent
$\sum_{j_0=1}^{\infty}2^{-j_0}m'_l(\w 2^{-j_0})
\prod_{j=1, j\neq j_0}^{\infty} m_l(\w 2^{-j}).$
\end{lem}

\textbf{Proof.}
Using the introduced notation we have
$$
\left|\left(\prod_{j=1}^{\infty}m_l\left(\frac{\w}{ 2^{j}}\right)\right)'-\left(\prod_{j=1}^{n}m_l\left(\frac{\w}{ 2^{j}}\right)\right)'\right|=
$$
$$
=
\left|\sum_{j_0=1}^{\infty}2^{-j_0}m'_l\left(\frac{\w}{ 2^{j_0}}\right)
\prod_{j=1, j\neq j_0}^{\infty} m_l\left(\frac{\w}{ 2^{j}}\right)-\sum_{j_0=1}^{n}2^{-j_0}m'_l\left(\frac{\w}{ 2^{j_0}}\right)
\prod_{j=1, j\neq j_0}^{n} m_l\left(\frac{\w}{ 2^{j}}\right)\right|\leq
$$
$$
\leq 
\left|
\sum_{j_0=n+1}^{\infty}2^{-j_0} m'_l\left(\frac{\w}{ 2^{j_0}}\right)
\prod_{j=1, j\neq j_0}^{\infty} m_l\left(\frac{\w}{ 2^{j}}\right)
\right|+
$$
$$
+
\left|
\sum_{j_0=1}^{n}2^{-j_0} m'_l\left(\frac{\w}{ 2^{j_0}}\right)\left(
\prod_{j=1, j\neq j_0}^{\infty} m_l\left(\frac{\w}{ 2^{j}}\right)-\prod_{j=1, j\neq j_0}^{n} m_l\left(\frac{\w}{ 2^{j}}\right)\right)
\right|=:I_{5,n}(\w)+I_{6,n}(\w).
$$
The application of the Lagrange Theorem and Lemma \ref{Cm'l} yields
$|m_l(\w)|\leq 1+A |\w|,$ where $A$ is a constant. Hence
$$
\left|\prod_{j=1}^{j_0-1} m_l\left(\frac{\w}{ 2^{j}}\right)\right|\leq
\prod_{j=1}^{j_0-1}\left(1+A|\w| 2^{-j}\right)=
e^{\sum_{j=1}^{j_0-1}\ln (2^j+A|\w|)-\ln 2^j}\leq
e^{A|\w|\sum_{j=1}^{j_0-1}\frac{1}{2^{j}}}\leq e^{A|\w|}.
$$
So using additionally Lemmas \ref{Cm'l} and \ref{cphiC} for the first sum we get
$$
I_{5,n}(\w)\leq
\sum_{j_0=n+1}^{\infty}2^{-j_0} \left|m'_l\left(\frac{\w}{ 2^{j}}\right)
\prod_{j=1}^{j_0-1} m_l\left(\frac{\w}{ 2^{j}}\right)\widehat{\varphi_l}\left(\frac{\w}{2^{j_0}}\right)\right|\leq
$$
$$
\leq  \left(\left\|(m^M)'\right\|_{C}+\left\|m_l'-(m^M)'\right\|_{C}\right) e^{A|\w|}
\left(\left\|\widehat{\varphi^M}\right\|_{C(\mathbb{R})}+
\left\|\widehat{\varphi^M}-\widehat{\varphi_l}\right\|_{C(\mathbb{R})}\right) 
2^{-n},
$$
where all factors are bounded as $a \leq \w \leq b,$ $l \in \mathbb{N}.$
Thus $I_{5,n}(\w)\rightarrow 0$ as $n \to \infty.$

For the second sum $I_{6,n}(\w)$ we obtain
$$
I_{6,n}(\w)= \left|\sum_{j_0=1}^{n}2^{-j_0} m'_l\left(\frac{\w}{ 2^{j_0}}\right) 
\prod_{j=1, j\neq j_0}^{n} m_l\left(\frac{\w}{ 2^{j}}\right)\left(\widehat{\varphi_l}\left(\frac{\w}{2^{n}}\right)-1\right)\right|.
$$
Since the function $\widehat{\varphi_l}$ is continuous, $\widehat{\varphi_l}(0)=1,$
and $a\leq |\w|\leq b,$ it follows that 
$\left|\widehat{\varphi_l}\left(\frac{\w}{2^{n}}\right)-1\right|=\ve_1(n),$  
$\ve_1(n)\rightarrow 0$ as $n \to \infty.$ 
So we get
$$
I_{6,n}(\w)\leq \left(1-\frac{1}{2^{n+1}}\right)
\left(\left\|(m^M)'\right\|_{C}+\left\|m_l'-(m^M)'\right\|_{C}\right) e^{A|\w|}
\ve_1(n),
$$
where all factors are bounded as $a \leq \w \leq b,$ $l \in \mathbb{N}.$
Thus $I_{5,n}(\w)\rightarrow 0$ as $n \to \infty$ $\Box.$

\section{Convergence of time and frequency radii for the wavelet functions}
%The equality  (\ref{qsw}) shows that $\widehat{\psi_l^{\bot}}$ is even. The function 
%$\widehat{\psi^M}$ is also even. Therefore $\w_{0\widehat{\psi_l^{\bot}}}=
%\w_{0\widehat{\psi^M}}=0$ 
%and $t_{0\psi_l^{\bot}}=
%t_{0\psi^M}=1/2.$
%Using the structure of formula (\ref{qsw}), applying Lemmas \ref{Cml}, \ref{cPhi}, \ref{cfr} % for the case of the frequency radii and applying additionally 
%Lemmas \ref{Cm'l}, \ref{cPhi'}, \ref{cphi'L}, and Remark \ref{cphi'CR}
%for the case of the time radii
%one can prove the convergence of the frequency and time radii for the wavelet functions in an %absolutely similar manner as Lemma \ref{cfr} and Lemma \ref{cti}. 
%It is only need to add the inequality estimating a difference of products, namely
%$|a_1 b_1-a_2b_2|\leq |a_1||b_1-b_2|+|b_2||a_1-a_2|$ and to apply it to the terms
%$
%\overline{m_l^{\bot}}\left(\w/2+\pi\right)
%	\widehat{\varphi^{\bot}_l}\left(\w/2\right)
%$ and
%$
%\overline{m^M}\left(\w/2+\pi\right)
%	\widehat{\varphi^M}\left(\w/2\right).
%$
%So we have
\begin{teo}\label{cfrtiw}
$
|\Delta^2_{\widehat{\psi_l^{\bot}}}-\Delta^2_{\widehat{\psi^M}}|=
O\left(\max\{\mu(l),\, l C_0^{-l+2\log_2 \frac {1+\ve(l)}{c}}\}\right),
$

\noindent
$
|\Delta^2_{\psi_l^{\bot}}-\Delta^2_{\psi^M}|=
O\left(\max\{\mu(l),\, (4 e^{2\w_0})^{-l+2\log_2 \frac {1+\ve(l)}{c}}\}\right)
$
as $l \to \infty.$
Parameters 
are defined by (\ref{param}).
\end{teo}
\textbf{Proof.}
The equality  (\ref{qsw}) shows that $\widehat{\psi_l^{\bot}}$ is even. The function 
$\widehat{\psi^M}$ is also even. Therefore $\w_{0\widehat{\psi_l^{\bot}}}=
\w_{0\widehat{\psi^M}}=0$ 
and $t_{0\psi_l^{\bot}}=
t_{0\psi^M}=1/2.$ The mask $m^{\bot}_l$ is real-valued function.
Therefore, using the stucture of formula (\ref{qsw}) and applying Lemmas \ref{Cml}, \ref{cPhi}, and Theorem \ref{cfr} we get for the frequency radii
$$
 |\Delta^2_{\widehat{\psi_l^{\bot}}}-\Delta^2_{\widehat{\psi^M}}|=
 $$
 $$
 =
\left|\int_{\mathbb{R}}\w^2\left((m^{\bot}_l)^2\left(\frac{\w}{2}+\pi\right)\left(\widehat{\varphi_l^{\bot}}\right)^2\left(\frac{\w}{2}\right)
-(m^M)^2\left(\frac{\w}{2}+\pi\right)\left(\widehat{\varphi^M}\right)^2\left(\frac{\w}{2}\right)\right)\, d\w\right|\leq
$$
$$
\leq
\|(m^{\bot}_l)^2\|_C
\int_{\mathbb{R}}\w^2\left|\left(\widehat{\varphi_l^{\bot}}\right)^2\left(\frac{\w}{2}\right)
-\left(\widehat{\varphi^M}\right)^2\left(\frac{\w}{2}\right)\right|\, d\w+
$$
$$
+\|(m^{\bot}_l)^2-(m^M)^2\|_C\int_{\mathbb{R}}\w^2\left(\widehat{\varphi^M}\right)^2\left(\frac{\w}{2}\right)\, d\w=
$$
$$
=O\left(|\Delta^2_{\widehat{\varphi_l^{\bot}}}-\Delta^2_{\widehat{\varphi^M}}|\right)+
O\left(\|\Phi_l-1\|_C\right)+
O\left(\|m_l-m^M\|_C\right)=
$$
$$
=O\left(\max\{\mu(l),\, l C_0^{-l+2\log_2 \frac {1+\ve(l)}{c}}\}\right).
$$
Going on to the time radii we use the identity
$\Delta^2_f=\int_{\mathbb{R}}t^2|f(t)|^2 \, dt-t^2_{0f}.$ 
Then applying  
Lemmas \ref{Cm'l}, \ref{cPhi'}, \ref{cphi'L} and Remark \ref{cphi'CR}
we obtain
$$
2 \pi 
|\Delta^2_{\psi_l^{\bot}}-\Delta^2_{\psi^M}|=
\left|\int_{\mathbb{R}}\left(\widehat{\psi_l^{\bot}}'\right)^2(\w)-
\left(\widehat{\psi^M}'\right)^2(\w) \, d\w\right|\leq
$$
$$
\leq \left\|\widehat{\psi_l^{\bot}}'+\widehat{\psi^M}'\right\|_{C(\mathbb{R})}
\int_{\mathbb{R}}\left|\widehat{\psi_l^{\bot}}'(\w)-
\widehat{\psi^M}'(\w)\right| \, d\w\leq
$$
$$
\leq
A \int_{\mathbb{R}}\left|
\left.m^{\bot}_l\right.'\widehat{\varphi_l^{\bot}} +
m^{\bot}_l\widehat{\varphi_l^{\bot}}' 
-i m^{\bot}_l\widehat{\varphi_l^{\bot}}
-
\left.m^M\right.'\widehat{\varphi^M} -
m^M\widehat{\varphi^M}'+i m^M\widehat{\varphi^M}\right|\leq
$$
$$
\leq A\left(
\|\left.m^{\bot}_l\right.'-i m^{\bot}_l\|_C \int_{\mathbb{R}}\left|\widehat{\varphi_l^{\bot}}-
\widehat{\varphi^M}\right| +\|\left.m^{\bot}_l\right.'-\left.m^M\right.'\|_{C}
\int_{\mathbb{R}}\left|\widehat{\varphi^M}\right|+\right.
$$
$$
\left.+
\|m^{\bot}_l\|_C \int_{\mathbb{R}}\left|\widehat{\varphi_l^{\bot}}'-
\widehat{\varphi^M}'\right| +\|m^{\bot}_l-m^M\|_{C}
\int_{\mathbb{R}}\left|\widehat{\left.\varphi^M\right.'}- i \widehat{\varphi^M}\right|
\right)=:I_{7,l}+I_{8,l}+I_{9,l}+I_{10,l}.
$$
The application of Lemmas \ref{cphiC}, \ref{cphiL}, and \ref{cPhi}
yields
$$
I_{7,l}=O\left(\|\Phi_l-1\|_{C}\right)+O\left((4e^{2\w_0})^{-l+2\log_2 \frac {1+\ve(l)}{c}}\right)=
O\left(\max\{\mu(l),\, (4 e^{2\w_0})^{-l+2\log_2 \frac {1+\ve(l)}{c}}\}\right).
$$
Using the Definition of $m^{\bot}_l,$ Lemmas \ref{cPhi'}, \ref{cPhi}, and \ref{Cm'l} we get
$$
I_{8,l}=O\left(\|\Phi'_l\|_C\right)+O\left(\|\Phi_l-1\|_C\right)+
O\left(\|m'_l-\left.m^M\right.'\|_C\right)=
$$
$$
=
O\left(\max\{\mu(l),\, (4 e^{2\w_0})^{-2l+4\log_2 \frac {1+\ve(l)}{c}}\}\right).
$$  
From Theorem \ref{cti} it follows that
$$
I_{9,l}=O\left(\max\{\mu(l),\, l C_0^{-l+2\log_2 \frac {1+\ve(l)}{c}}\}\right).
$$
Finally, using Lemmas \ref{cPhi} and \ref{Cml} we get
$$
I_{10,l}=O\left(\|\Phi_l-1\|_C\right)+O\left(\|m_l-m^M\|_C\right)=
O\left(\max\{\mu(l),\, (4 e^{2\w_0})^{-2l+4\log_2 \frac {1+\ve(l)}{c}}\}\right).
$$
Thus, collecting the estimations we obtain
$$
|\Delta^2_{\psi_l^{\bot}}-\Delta^2_{\psi^M}|=
O\left(\max\{\mu(l),\, (4 e^{2\w_0})^{-l+2\log_2 \frac {1+\ve(l)}{c}}\}\right) \Box.
$$

%\section{Conclusion.}
%So, the family of quasispline wavelets   is constructed.
%A non-orthogonal quasispline mask is defined in  (\ref{ml}). The construction use  
%a linear method of summation satisfying the conditions (\ref{con1})-(\ref{con3}).
%The scaling and wavelet functions  exponentially decay (see Theorem \ref{expp}) and theirs Fourier transforms have %decay of polynomial order  (see Theorem \ref{Smo}).
%The squares of the uncertainty constants of the scaling (wavelet) functions tend to those of the Meyer scaling %(wavelet)  function
%(see Theorems \ref{cfr}, \ref{cti}, \ref{cfrtiw}). So, the quasispline wavelets solve the problem (1).

\section*{Acknowledgments.}
The author thanks Professor I. Ya. Novikov for setting the problem and Professor M. A. Skopina for valuable discussions.




\begin{thebibliography}{00}


\bibitem{ChW}
Chui C.K., Wang J. High-order orthonormal scaling functions and wavelets give poor time-frequency localization // The J. of Fourier Anal. and Appl., Vol. 2, No. 5, 1996, pp. 415-426.
 
 \bibitem{ChW2}
  Chui C.K., Wang. J. A study of asymptotically optimal time-frequency
localization by scaling functions and wavelets // Annals of Numerical Mathematics, Vol.
4, 1997, pp. 193-216.

\bibitem{GL}
Goodman T.N.T., Lee S. L., Asymptotic Optimality in Time-Frequency
Localization of Scaling Functions and Wavelets.  Frontiers in Interpolation and
Application, (Eds.) N.K. Govil, H.N. Mhaskar, R.N. Mohapatra, Z. Nashed, and
J. Szadados, 2006.
 
 
 \bibitem{N1}
 Novikov I. Ya. Modified Daubechies wavelets preserving loca\-lization with growth of smoothness // East J. Approximation,  Vol. 1, No. 3, 1995, pp. 314-348.
  
  \bibitem{N2}
 Novikov I.  Ya. Uncertainly constants for modified Daubechies wavelets [in Russian] // Izv. Tul'sk. Gos. Univ. Ser. Mat. Mekh. Inform., Vol. 4, No. 1, 1998, pp. 107-111.
  
  \bibitem{LSMZ}
 Lebedeva E.A. Exponentially decaying wavelets with uncertainty constants uniformly bounded with respect to the smoothness para\-meter // Sib. Math. J., Vol. 49, No. 3, 2008, pp. 457-473.
 
 \bibitem{D} Daubechies I. Ten lectures on wavelets. CBMS-NSF Series in Appl. Math., SIAM, Philadelphia, 1992.
 
 \bibitem{NPS}
 Novikov I. Ya., Protasov V. Yu., and Skopina M. A. The Theory of Wavelets [in Russian]. Fizmathlit, Moscow (2005).
 
\end{thebibliography}
\end{document}